\documentclass[12pt,reqno]{amsproc}
\allowdisplaybreaks
\numberwithin{equation}{section}

\usepackage[hyphens]{url}
\makeatletter
\g@addto@macro{\UrlBreaks}{\do\/\do\-}
\makeatother
\usepackage{empheq}
\usepackage{amssymb}
\usepackage{stmaryrd}

\usepackage{enumerate}
\usepackage{fullpage}
\usepackage{textcomp}
\usepackage{lettrine}

\usepackage{graphicx}
\usepackage{subfigure}

\usepackage{morefloats}
\usepackage[font=small,labelfont=bf,
   justification=justified,
   format=plain]{caption} % 'format=plain' avoids hanging indentation
\usepackage[debug=false, colorlinks=true, pdfstartview=FitV,
linkcolor=blue, citecolor=blue, urlcolor=magenta,breaklinks=true]{hyperref}
\usepackage[semicolon,square,authoryear,sort]{natbib}

\usepackage[doipre={DOI:~}]{uri}

\usepackage{color}
\usepackage{colortbl}
\usepackage[most]{tcolorbox}

\newtcbox{\mymath}[1][]{%
    nobeforeafter, math upper, tcbox raise base,
    enhanced, colframe=blue!30!black,
    colback=blue!30, boxrule=1pt,
    #1}

\newlength{\drop}
\definecolor{amethyst}{rgb}{0.6, 0.4, 0.8}
\definecolor{burgundy}{rgb}{0.5, 0.0, 0.13}

\usepackage{longtable} 
\usepackage{multirow} 
\usepackage{rotating} 
\usepackage{bigstrut} 
\usepackage{hhline} 
\usepackage{amsthm}

\newtheoremstyle{remboldstyle}
  {}{}{}{}{\bfseries}{.}{.5em}{{\thmname{#1 }}{\thmnumber{#2}}{\thmnote{ (#3)}}}
\theoremstyle{remboldstyle}

\title{\textbf{Reactive Transport Modeling with Physics-Informed Machine Learning for Critical Minerals Applications}}

\author{\textbf{Kripa~Adhikari} \\
{\small Department of Civil and Environmental Engineering \\
University of Houston, Houston, Texas 77204, USA.}\\
\textbf{Md.~Lal Mamud}, \textbf{Maruti Kumar Mudunuru} \\
{\small Pacific Northwest National Laboratory, Richland, Washington 99352, USA.}\\
\textbf{Kalyan~B.~Nakshatrala}\\
{\small Department of Civil and Environmental Engineering \\
University of Houston, Houston, Texas 77204, USA.}\\
{\small Correspondence to: \texttt{knakshatrala@uh.edu} (email), +1-713-743-4418 (phone); \texttt{maruti@pnnl.gov} (email), +1-509-375-6645 (phone)}}

\keywords{physics-informed neural networks; 
tailored weighting; 
reactive-transport; 
flow through porous media;
patch tests;
critical mineral science}

\begin{document}

%===========================;
%  Title page of the paper  ;
%===========================;
\begin{titlepage}
  \drop=0.1\textheight
  \centering
  \vspace*{\baselineskip}
  \rule{\textwidth}{1.6pt}\vspace*{-\baselineskip}\vspace*{1pt}
  \rule{\textwidth}{0.4pt}\\[\baselineskip]
       {\Large \textbf{\color{burgundy}
       Reactive Transport Modeling with Physics-Informed Machine Learning for Critical Minerals Applications}}\\[0.1\baselineskip]
       \rule{\textwidth}{0.4pt}\vspace*{-\baselineskip}\vspace{1pt}
       \rule{\textwidth}{1.6pt}\\[\baselineskip]
       \scshape
       %
       % An e-print of this paper is available on arXiv. \par
       %
       \vspace*{-1\baselineskip}
       Authored by \\[-0.2\baselineskip]

   {\large K.~Adhikari\par}
  {\itshape Graduate Student, Department of Civil \& Environmental Engineering, \\
  University of Houston, Houston, Texas 77204.}\\[0\baselineskip]

  {\large Md.~L.~Mamud\par}
  {\itshape Postdoctoral Associate, Pacific Northwest National Laboratory.}\\[0\baselineskip]

  {\large M.~K.~Mudunuru\par}
  {\itshape Staff Scientist, Pacific Northwest National Laboratory.}\\[0\baselineskip]

  {\large K.~B.~Nakshatrala\par}
  {\itshape Department of Civil \& Environmental Engineering \\
  University of Houston, Houston, Texas 77204.} 

%=====================;
%  Front page figure  ;
%=====================;
\begin{figure*}[ht]
    \centering
    \includegraphics[scale=0.55]{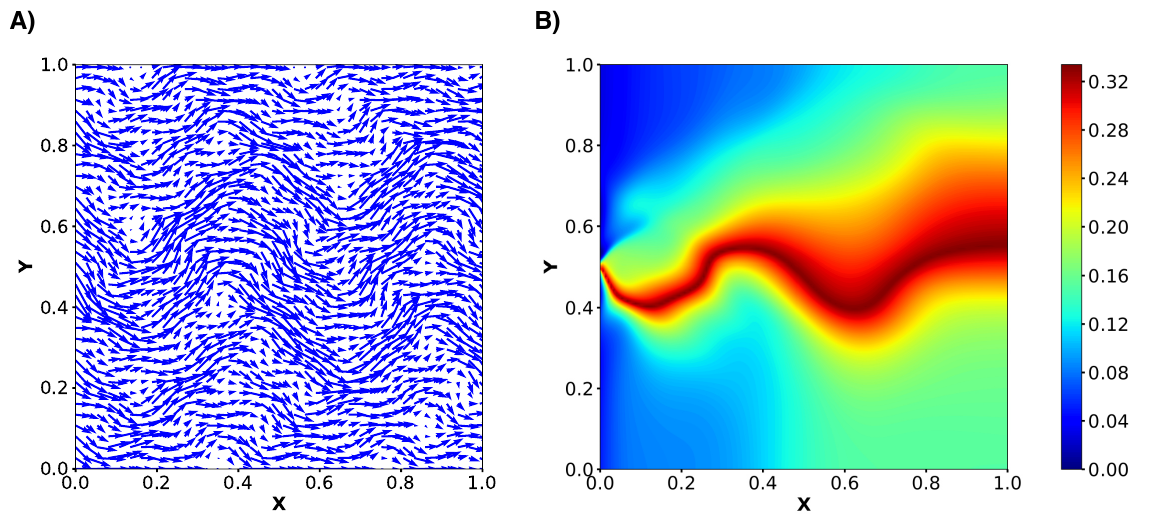}
    \vspace{-0.2in}   
    \captionsetup{labelformat=empty}
    \caption{{Figure \small (A) shows the velocity field, with arrows indicating variations in direction and magnitude typical of an \emph{in situ} leaching scenario. Flow through heterogeneous media creates channeling and recirculation zones that strongly affect reagent transport and mixing. (B) depicts the plume of product C (e.g., a metal–ligand complex) formed by a fast bimolecular reaction between reactants A (e.g., acid donor) and B (e.g., complexing agent). Concentrations are predicted using physics-informed neural networks (PINNs). The plume starts at the left boundary, where reactants enter, and sharpens moving right. Capturing these mixing-limited fronts is key to optimizing reagent injection, maximizing critical mineral extraction/recovery, and minimizing reagent use and by-products.}}
\end{figure*}
\vspace{-0.22in}
%  \vfill
  {\scshape 2025} \\
  {\small Computational \& Applied Mechanics Laboratory} \par
\end{titlepage}

%=========================;
%  Abstract of the paper  ;
%=========================;
\begin{abstract} 
    This study presents a physics-informed neural networks (PINNs) framework for reactive transport modeling for simulating fast bimolecular reactions in porous media. Accurate characterization of chemical interactions and product formation in surface and subsurface environments is essential for advancing critical mineral extraction and related geoscience applications. The proposed methodology sequentially addresses the flow and diffusion–reaction subproblems. The flow field is computed using a mixed formulation, while the diffusion–reaction system is modeled via two uncoupled tensorial diffusion equations reformulated in terms of chemical invariants. PINNs are employed to solve the governing equations, enabling data-efficient, mesh-free prediction of chemical concentration fields. The framework is validated through a series of benchmark problems involving flow in heterogeneous porous media. Initial verification is conducted using patch tests for the flow field, followed by validation of the transport problem with emphasis on preserving non-negativity of concentrations. The complete fast bimolecular reaction scenario is then solved, yielding spatial distributions of reactants and product species. Results demonstrate that the PINNs-based approach effectively captures sharp, mixing-limited reaction fronts and dispersive mixing behavior, offering reliable predictions of reactive plume evolution. These capabilities are crucial for evaluating long-term subsurface behavior in applications such as fluid storage, energy extraction, and efficient extraction of critical minerals.
\end{abstract}

\maketitle

%==================================;
%  Include all the sections below  ;
%==================================;
\setcounter{figure}{0}   

%*********************************************;
%                                             ;
%  NAME                                       ;
%    S1_PINNS_Fast_Intro.tex                  ;
%                                             ;
%*********************************************;
\section{INTRODUCTION AND MOTIVATION}
\label{Sec:S1_PINNS_Fast_Intro}
\begin{figure}[hbt!]
    \centering
    \includegraphics[scale = 0.75]{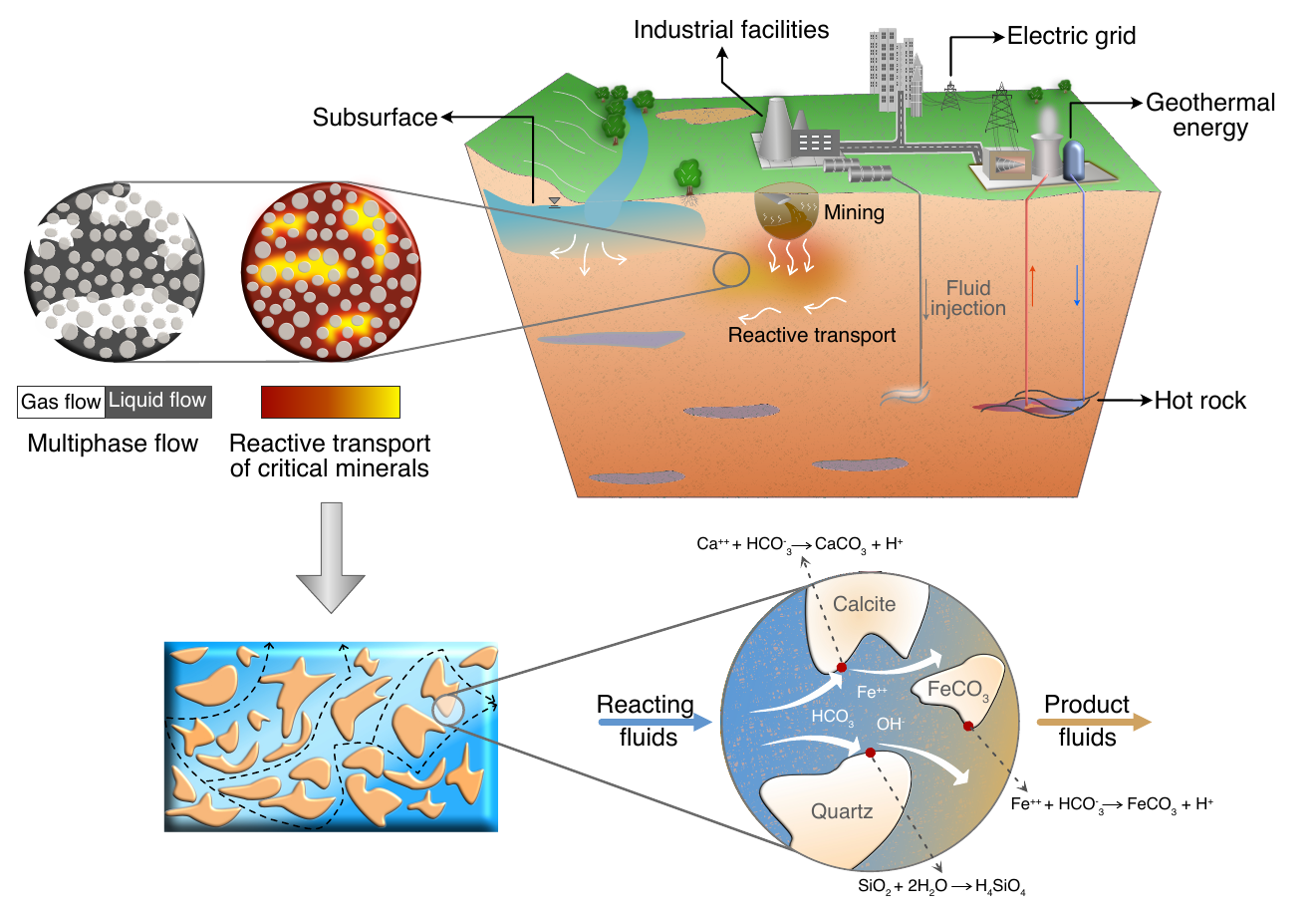}
    \caption{\textsf{Subsurface mass and reactive-transport for critical minerals and materials:} The figure illustrates a range of scenarios involving fluid (gas and liquid) flow and transport within the Earth's subsurface---one of the largest porous media for fluid flow. The figure showcases various subsurface processes, including groundwater movement, renewable energy systems such as geothermal energy, geologic carbon storage, and reactive transport. The mineral phases can be clays, carbonates, Fe/Mn oxides, phosphates, and silicates. For instance within the critical mineral context \citep{cantrell1987rare}, this may include $\text{Clay--Na} + \text{Li}^+ \rightarrow \text{Clay--Li} + \text{Na}^+$ or $\text{REE}^{3+} + \text{CO}_3^{2-} \rightarrow \text{REECO}_3^{+}$ within a much bigger multi-step reaction settings. This reaction focuses on a bimolecular complexation step in which a rare earth element (REE) forms a 1:1 carbonate complex in aqueous solution, commonly observed in alkaline or carbonate-rich environments.}
    \label{fig:mass_transport}
\end{figure}

\lettrine[findent=2pt]{\fbox{\textbf{T}}}{he} mixing of chemicals and the formation of products is important in both surface and subsurface environments. At the surface, chemical mixing affects water resources \citep{tripathi2021contamination}, while underground, it influences processes such as fossil fuel recovery, mineral extraction, and geothermal energy through interactions with subsurface chemicals \citep{national2005contaminants}. Many subsurface processes involve mass transport, as the subsurface acts as a vast porous medium that facilitates fluid flow. These processes encompass not only groundwater resources but also reactive transport due to mineral-water-rock interactions, fluid flow in renewable energy systems like geothermal energy, and geologic carbon storage, as illustrated in Fig.~\ref{fig:mass_transport}. 

Understanding mass transport is crucial for predicting plume behavior when flow, diffusion, and reactions occur. In many cases, chemical reactions are sufficiently fast (e.g., \emph{in situ} leaching scenarios), allowing them to be approximated as instantaneous \citep{dentz2011mixing}. Mixing-limited reactions in porous media arise when reaction kinetics are fast relative to transport, such that reaction rates are controlled by incomplete mixing driven by diffusion and dispersion. This behavior has been extensively studied using theoretical, numerical, and experimental approaches, highlighting the strong dependence of reactive plume evolution on pore-scale heterogeneity and mixing dynamics rather than intrinsic reaction rates. Previous studies have employed pore-scale simulations, upscaled continuum models, and particle- and volume-averaged approaches to quantify mixing-controlled reaction rates and effective reaction terms \citep{valocchi2019mixing}. However, accurately capturing sharp reaction fronts and dispersive mixing at the continuum scale remains challenging due to numerical diffusion, scale dependence, and the need for fine spatial resolution. Comprehensive reviews of these challenges and modeling strategies are provided by. In this study, flow through porous media, chemical transport, and fast bimolecular reactions are treated as interconnected processes: the flow field informs transport and reaction dynamics, and all three components together determine reactive transport behavior. The modeling framework proceeds in three coupled stages: first solving the flow problem to obtain the velocity field, then computing the dispersion tensor for transport, and finally solving the steady-state diffusion-reaction equation to model fast chemical reactions.

Darcy equations are widely used to model Newtonian fluid flow in porous media across a range of applications and naturally incorporate anisotropy through a tensorial permeability field \citep{neuman1977theoretical}. There are two broad approaches to solving Darcy equations. The first, the single-field formulation, rewrites the governing equations solely in terms of pressure. This approach has two main drawbacks: (a) many generalizations of Darcy equations (e.g., the Darcy–Brinkman and Darcy–Forchheimer models) cannot be expressed solely in terms of pressure \citep{chang2017modification}; and (b) the accuracy of the velocity is typically poor, as it is derived from the pressure field, even though the velocity is often more important than the pressure. The second approach, mixed-field formulations, discretizes both pressure and velocity simultaneously, maintaining good accuracy in both solution fields \citep{masud2002stabilized,nakshatrala2006stabilized,arbogast2007multiscale, nakshatrala2011numerical}. Due to these advantages, mixed formulations have been widely applied in flow through porous media applications, including subsurface energy systems and reservoir simulations \citep{chavent1984discontinuous, ewing1984mixed}. Therefore, this study adopts a mixed formulation, as we want to predict accurate velocity fields, especially in heterogeneous porous media, where capturing flux continuity is critical. As we demonstrate later, using a mixed formulation is quintessential for the proposed modeling framework to satisfy patch tests---simple canonical problems commonly used to assess the accuracy of numerical methods. 

%--------------------------------;
%  Figure 1: Modeling evolution  ;
%--------------------------------;
\begin{figure}[hbt!]
    \centering
    \includegraphics[scale = 0.7]{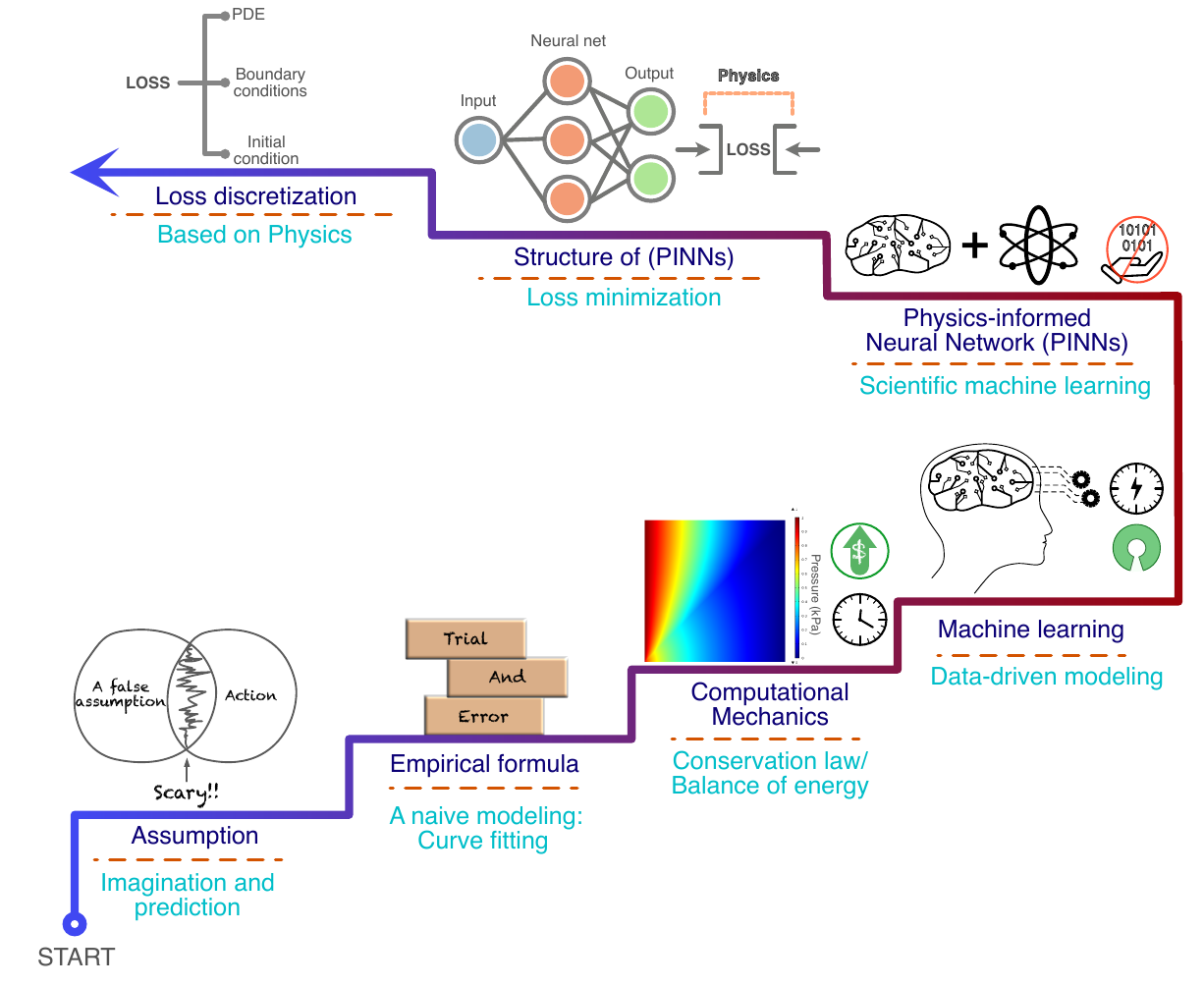}
    \caption{\textsf{Modeling evolution:} The figure illustrates the progression of modeling approaches over time. It begins with basic assumptions and advances to curve fitting. With the deeper understanding of physics and mechanics, simulation and computational mechanics emerged. Eventually, with the advent of advanced programming languages, data-driven and physics-informed machine learning has experienced significant growth.}
    \label{fig:pinns}
\end{figure}

As illustrated in Fig.~\ref{fig:pinns}, modeling has evolved through several stages: from naïve guesswork to empirical formulations, and later to computational mechanics, which enabled more precise solutions across applications. Although modern programming languages and open-source libraries have improved accessibility, developing robust, scalable solvers for complex multiphysics problems still demands substantial computational resources. This has spurred interest in machine learning, whose data-driven approaches excel when large datasets are available, ranging from image recognition \citep{krizhevsky2017imagenet} to human sciences \citep{lake2015human}, driven by increased computational power and improved algorithms \citep{jordan2015machine}. Yet, many fields face limited or unavailable data \citep{latrach2024critical}, and solving nonlinear partial differental equations (PDEs) for coupled physics problems remains challenging. Physics-informed neural networks (PINNs) address this by embedding governing physical laws into the training loss, enabling modeling and prediction even when data are scarce or incomplete \citep{raissi2018hidden, raissi2019physics, zhang2019quantifying, raissi2019deep, yang2019adversarial}.

In porous media flow, various modeling techniques have been developed to tackle complex problems. Widely used methods include the finite element \citep{diersch2013feflow, fang2017coupled, pinder2006subsurface}, finite volume \citep{terekhov2022finite, burbulla2022finite, versteeg2007introduction}, and finite difference methods \citep{zhan2023generalized, banaei2021numerical, smith1985numerical}. With advances in computational programming, machine learning (ML) approaches—especially after 2015 \citep{delpisheh2024leveraging}—have gained traction for predicting flow fields, estimating recovery efficiency, and learning effective constitutive relationships. Architectures such as artificial neural networks (ANNs), convolutional neural networks (CNNs), and other supervised methods have been explored depending on the problem and data availability \citep{luo2023review, hanga2019machine, hegde2017use}. However, conventional ML is limited by the scarcity of training data in porous media applications.

Recent studies have used machine learning to accelerate reactive transport simulations by replacing chemical solvers with trained surrogates or reduced models \citep{leal2020accelerating, kyas2022accelerated}. While effective within operator-splitting frameworks, these approaches treat transport and reactions separately. In contrast, PINNs approach used in this work integrates flow, transport, and reactions in a coupled framework. Hence, PINNs can resolve governing physics that are enforced directly within each subproblem using PINNs, and the overall reactive transport process is solved in a sequential, physics-consistent manner, where the flow solution provides inputs to the reactive transport model and provide a flexible and well-suited framework for inverse modeling and data-scarce scenarios. Physics-informed machine learning (PIML) has been applied to porous media flow \citep{fraces2021physics, mao2020physics, wang2020deep}, from subsurface energy systems \citep{latrach2024critical} to vascular cooling \citep{jagtap2023coolpinns}. Though termed PIML or PINNs in the literature \citep{harp2021feasibility, jagtap2023coolpinns}, both share the same principles; we adopt ``PINNs'' throughout this paper.

Researchers are increasingly adopting PINNs across diverse applications, including fluid dynamics \citep{cai2021physics}, reservoir pressure management \citep{harp2021feasibility}, subsurface energy systems \citep{shima2024modeling, latrach2024critical}, and biomedical flows \citep{yin2021non}. Despite their potential, PINNs face challenges in modeling nonlinear two-phase transport in porous media, often yielding inaccurate approximations \citep{fuks2020limitations}. While convolutional neural networks (CNNs) have effectively solved problems such as elliptic pressure equations with discontinuous coefficients \citep{makauskas2023comparative}, the use of PINNs for fast bimolecular reactions with fully coupled flow-transport-reaction mechanisms remains unexplored. Previous studies have addressed diffusion-reaction problems using finite element methods \citep{mudunuru2016enforcing} or data-driven ML \citep{ahmmed2021comparative, vesselinov2019unsupervised}, but PINNs have yet to be applied in this context.

Addressing these gaps, this study aims to extend the application of PINNs and assess their suitability for reactive-transport modeling in the context of critical mineral science. To this end, we focus on three progressively complex problems: First, we consider a heterogeneous flow subproblem, featuring three patch tests---horizontal, vertical, and inclined---solved using PINNs and validated against finite element method (FEM) solutions. This addresses the question:
\begin{enumerate}
    \item[(Q1)] How do PINNs perform in modeling flow through porous media under heterogeneous conditions?
\end{enumerate}
Next, we examine a transport problem governed by the anisotropic diffusion equation. The non-negativity of the solution is evaluated and compared with FEM results, addressing the question:
\begin{enumerate}
    \item[(Q2)] Does the PINNs formulation obey mathematical qualitative properties, such as the maximum principle?
\end{enumerate}
Finally, we investigate fast bimolecular reaction problems under both uniform and non-uniform flow, exploring whether PINNs can accurately capture coupled flow-transport-reaction dynamics:
\begin{enumerate}
    \item[(Q3)] Can PINNs solve fast bimolecular reaction problems, incorporating both uniform and non-uniform flow velocity?
\end{enumerate}

To this end, a full flow–transport–reaction model is implemented for uniform flow, with an explicit non-uniform velocity introduced to evaluate the framework under more complex reaction–diffusion conditions. The scientific contributions of this paper—including the mathematical modeling, computational framework, and technical insights—are highly relevant to applications in critical mineral science (e.g., see the discussion in \S\ref{SubSec:S2_DiffRxn}). This work advances our ability to simulate and predict subsurface chemical interactions by integrating geochemistry, fluid dynamics, and machine learning, providing a powerful toolset for informed decision-making in mineral exploration and resource recovery.

The rest of this paper is organized as follows. We explain the PINNs framework and introduce a set of governing equations that describe firstly the flow problem, second the transport problem, and then the diffusion-reaction problem for fast bimolecular reaction (\S\ref{Sec:S2_PINNS_Fast_GE}). Next, we present the verification of the PINNs for the flow subproblem with different patch tests (\S\ref{Sec:S3_PINNS_Fast_flow}). Then, we utilize PINNs to solve for non-negative solutions under the transport problems (\S\ref{Sec:S4_PINNS_Fast_Transport}). Lastly, we leverage the PINNs to solve fast bimolecular reactions between chemicals for uniform flow velocity and a non-uniform explicit velocity (\S\ref{Sec:S5_PINNS_Fast_NR}). Finally, we draw conclusion remarks and propose potential future work (\S\ref{Sec:S6_PINNS_Fast_Closure}).                

%*********************************************;
%                                             ;
%  NAME                                       ;
%    S2_PNNS_Fast_GE.tex                      ;
%                                             ;
%*********************************************;
\section{MODELING FRAMEWORK AND GOVERNING EQUATIONS}
\label{Sec:S2_PINNS_Fast_GE}
\subsection{PINNs modeling framework}
\begin{figure}[h]
    \centering
    \includegraphics[scale = 0.7]{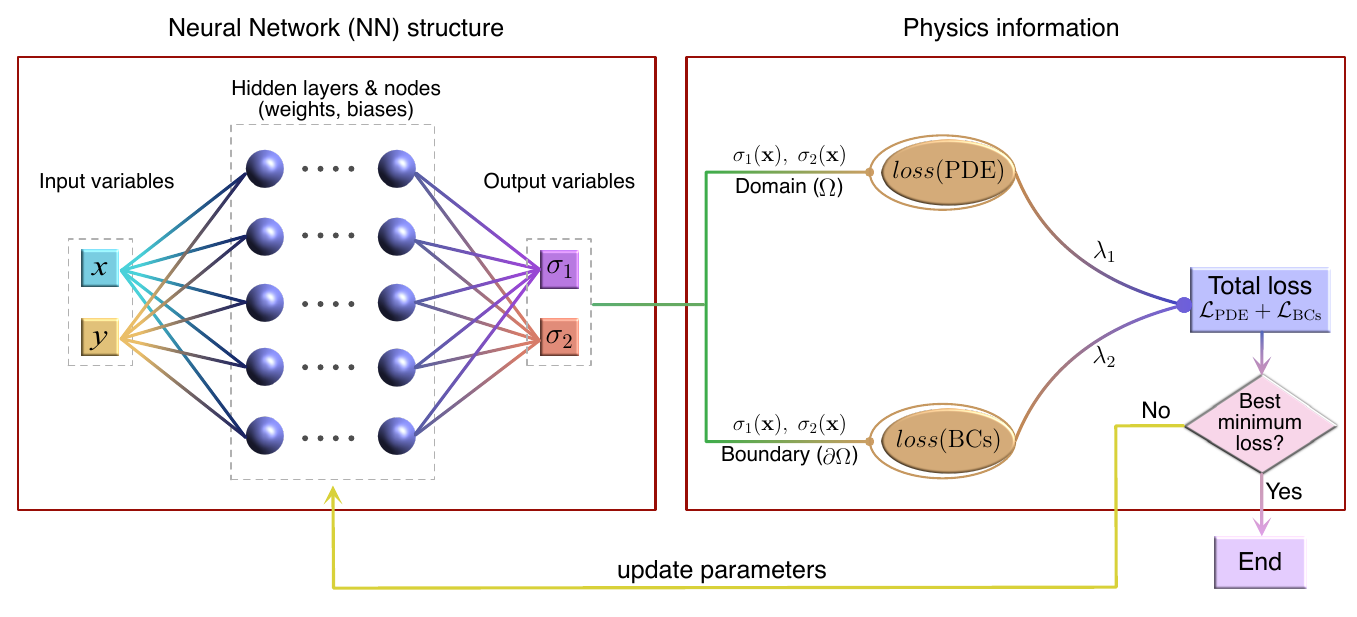}
    \caption{\textsf{PINNs framework:} The figure depicts the PINNs framework, which consists of a neural network (NN) structure and a physics-informed component. The input variables are processed through the NN to generate output variables. These outputs are evaluated against the governing equations of the system, including the PDE and boundary conditions. Residual losses are computed, and the weighted sum of these losses determines the total loss. If the total loss reaches its optimal minimum, the training concludes. Otherwise, the hyperparameters are updated based on the specified learning rate, and training continues until the maximum number of epochs is reached.}
    \label{fig:Fast_PINNs_framework}
\end{figure}

The fundamental principle of physics-informed neural networks (PINNs) is that the solution must satisfy the governing equations at each collocation point. Here, the domain is discretized into collocation points, and their coordinates are used as input variables for the model. The training process begins by initializing hyperparameters ($\theta$), including weights ($w$) and biases ($b$), using a random number generator. The neural network (NN) structure is then defined with a specified number of hidden layers and nodes as illustrated in Fig.~\ref{fig:Fast_PINNs_framework}. Since the architecture consists of more than two hidden layers and involves multiple layers of trainable weights, it can be referred to as a deep neural network (DNN) \citep{lecun2015deep}.

If $\mathbf{x}$ represents spatial points within the domain, let $n$ denote the total number of collocation points. The continuous coordinates in the domain are referred to as spatial points, while the discrete subset used to evaluate the PINNs residuals is referred to as collocation points. The input layer receives a total of $n$ collocation points, each with spatial coordinates 
($x$, $y$), forming an input tensor of shape ($n$, 2) for two dimensional (2D) problems. For each collocation point, a random value---generated using the default random value generator in the JAX library--- is assigned as the first step of the training process, which is then passed through the NN structure. At each neuron, an activation function is applied, and the output variables are computed as the initialized weights and biases propagate through the specified hidden layers ($l$) and neurons/nodes ($n_l$). Mathematically, the operation inside each neuron can be expressed as follows:
\begin{align}
    & a^l = \sum_{i=1}^{n_l} w_i^l \, \mathbf{x}_i^{(l-1)} + b^l\\
    &\sigma^l = f(a)^l
\end{align}
where, the spatial point $\mathbf{x}$ for a 2D domain consists of coordinates $x$ and $y$, $a$ represents the pre-activation output, and $f(\cdot)$ is the nonlinear activation function \citep{bishop2023deep}. Hyperbolic tangent ($\mathrm{tanh}$) is selected as the activation function in this study, given the smoothness and differentiability in solving the PDE within the PINNs framework \citep{raissi2019physics}. The term $\sigma^l$ denotes the output of the layer, also referred to as the activation. If $L$ represents the final output layer, then the output $\sigma_k$ represents the output of the whole neural network (NN) structure such that:
\begin{align}
    & a^L = \sum_{i=1}^{nL} w_i^L \, \mathbf{x}_i^{(L-1)} + b^L\\
    &\sigma_k = f(a)^L
\end{align}
where, $k \in \mathbb{N}$ ($\mathbb{N} = 1,2,3,...$), depends on the number of outputs, which varies based on the nature of the problem.
%----------------------------;
%  Figure 2: PINNa Framework ;
%----------------------------;

The resulting output variables then undergo physics evaluation inside the domain ($\Omega$) and at the boundary ($\partial\Omega$). Here, physics evaluation refers to the enforcement of physical laws to the entire domain to compute total residual loss. Within the domain, residual losses are computed based on the governing PDE, while boundary points are assessed against the prescribed boundary conditions, with the corresponding mean square error (MSE) formulas: 
\begin{align}
    & \mathcal{L}_{\mathrm{PDE}} = \frac{1}{n_\Omega} \; \sum_{i=1}^{n_\Omega}Residual \; \mathrm{PDE} \; (\mathbf{x}_i)^2 \\
    & \mathcal{L}_{\mathrm{BCs}} = \frac{1}{n_{\partial\Omega}} \; \sum_{j=1}^{n_{\partial\Omega}} \{true (\mathbf{x}^{\partial\Omega}_j) - predicted(\mathbf{x}^{\partial\Omega}_j)\}^2
\end{align}
where the $true (\mathbf{x}_j)$ denotes the prescribed Dirichlet or Neumann boundary condition values, specified during problem setup based on the nature of the PDE and boundary type. $n_\Omega$ are collocation points inside domain and $n_{\partial\Omega}$ are collocation points on the boundary. Each loss component is assigned a specific weight $\lambda_j$, where $j \in \mathbb{N}$, and depends on the number of loss terms. The total loss is then computed using the following expression:
\begin{align}
    \label{total_loss}
    \mathrm{Total \; loss} \; (\mathcal{L}) = \lambda_1 \cdot \mathcal{L}_\mathrm{PDE} + \lambda_2 \cdot \mathcal{L}_\mathrm{BCs} 
\end{align}
where $\lambda_1$ and $\lambda_2$ are the weights of residual loss due to the governing PDE $(\mathcal{L}_\mathrm{PDE})$ and boundary conditions $(\mathcal{L}_\mathrm{BCs})$, respectively. The weight values are crucial in this framework, as they determine the contribution of each loss term during the training process. They significantly influence the training outcome, determining whether the solution successfully satisfies all governing equations and boundary conditions \citep{radin2023effects}. There is no specific method to determine these loss weights in advance \citep{rao2021physics}, and they are typically determined through trial and error, as in this study. 

Once the total loss ($\mathcal{L}$) is computed, the training process evaluates whether it has reached the optimal minimum loss---ensuring it is lower than previous or upcoming losses. If the condition is satisfied, training stops; otherwise, the hyperparameters are iteratively updated using gradient-based optimization techniques \citep{bishop2023deep}---a well-suited approach for the smooth and differentiable loss functions. This involves efficiently computing the derivative of the error function with respect to the hyperparameters through backpropagation ($\nabla_\theta = \frac{\partial \mathcal{L}}{\partial \theta}$) \citep{bottou2012stochastic}. The gradient $\nabla_\theta$ determines the direction of optimization, while the step size, known as the learning rate ($\eta$), controls how far the update moves in that direction. The next value of the hyperparameters is computed using the following expression:
\begin{align}
    \theta_{t+1} = \theta_t - \eta \; \nabla_\theta
\end{align}
where, $t$ is the number of iterations. The learning rate $\eta$ plays a crucial role in optimizing the hyperparameters. An adaptive learning rate strategy is used with the Adam optimizer, as it offers a good tradeoff between training efficiency and robustness for the problems considered in this study \citep{wang2021understanding}. In contrast, a large learning rate may cause the updates to overshoot the optimum. The training process continues until either the best minimum loss is achieved or the maximum number of epochs is reached. An epoch refers to a single pass through the full set of training points during which the model updates its parameters to minimize the loss function. Once training is complete, the optimized parameters are used to compute the final output variables, which are then employed for visualization as required. In all subsequent studies, the same framework is utilized, with the primary variation being the governing equations and the specific output variables of interest.

%==========================================================;
%  Subsubsection: Flow through heterogeneous porous media  ;
%==========================================================;
\subsection{Flow through heterogeneous porous media}
To solve for a fast bimolecular reaction, it is necessary to first address the flow through the porous medium. A stabilized mixed formulation is employed to model the flow in porous media. We consider a 2D domain represented as $\Omega$ whose boundary is denoted by $\partial \Omega$ which is decomposed into complementary parts: $\Gamma^{p}$ and $\Gamma^{v}$. For  mathematical well-posedness, we require that $\Gamma^p \cup \Gamma^v = \partial \Omega$ and $\Gamma^p \cap \Gamma^v = \emptyset$. At any spatial point $\mathbf{x}$, the governing equations for the flow subproblem are given as follows:
%--------------------------------------;
%  Equation: Flow governing equations  ;
%--------------------------------------;
\begin{alignat}{2}
    \label{Eqn:PINNS_Fast_BoLM}
    &\mu \, \mathbf{K}^{-1}(\mathbf{x})
    \, \mathbf{v}(\mathbf{x}) 
    + \mathrm{grad}\big[p(\mathbf{x})\big] 
    = \rho \, \mathbf{b}(\mathbf{x}) 
    &&\quad \mathrm{in} \; \Omega \\
    \label{Eqn:PINNS_Fast_BoM}
    &\mathrm{div}\big[\mathbf{v}(\mathbf{x})\big] 
    = 0  
    &&\quad \mathrm{in} \; \Omega \\
    \label{Eqn:PINNS_pBC}
    &p(\mathbf{x}) = p^{\mathrm{p}}(\mathbf{x}) 
    &&\quad \mathrm{on} \; \Gamma^{p} \\
    \label{Eqn:PINNS_vBC}
    &\mathbf{v}(\mathbf{x}) \bullet 
    \widehat{\mathbf{n}}(\mathbf{x})  
    = v^{\mathrm{p}}(\mathbf{x}) 
    &&\quad \mathrm{on} \; \Gamma^{v} 
\end{alignat}
where $\rho$ is the fluid's density, $\mathbf{b}(\mathbf{x})$ the specific body force field, and $\mathbf{K}(\mathbf{x})$ is the coefficient of permeability. $\mu$ is the dynamic viscosity of the fluid and $p(\mathbf{x})$ is the pressure distribution at any spatial point $\mathbf{x}$. $p^{\mathrm{p}}(\mathbf{x})$ be the prescribed pressure on the boundary and $v^\mathrm{p}(\mathbf{x})$ 
 be the prescribed normal component of the velocity on the boundary where $\mathbf{v}(\mathbf{x})$ is the velocity field and $\widehat{\mathbf{n}}(\mathbf{x})$ is the unit normal outward. $\mathrm{grad}[\cdot]$ and $\mathrm{div}[\cdot]$ represent the spatial gradient and divergence operators, respectively. 

The first equation (i.e., Eq.~\eqref{Eqn:PINNS_Fast_BoLM}) represents the balance of linear momentum incorporating the Darcy model, Eq.~\eqref{Eqn:PINNS_Fast_BoM} accounts for mass conservation assuming incompressibility of the fluid, Eq.~\eqref{Eqn:PINNS_pBC} describes the pressure boundary condition, and the last equation (i.e., Eq.~\eqref{Eqn:PINNS_vBC}) specifies the velocity boundary condition, with the normal component of the velocity being prescribed.

%====================================;
%  Subsubsection: Transport problem  ;
%====================================;
\subsection{Transport problem}
Transport processes in porous media are central to subsurface applications, involving the movement of chemical species through a combination of advection and diffusion. Given the inherent heterogeneity in natural porous media and the possible anisotropy in diffusion, accurately capturing chemical concentration distributions is critical for understanding and predicting subsurface behavior. In this study, we employ a classical anisotropic diffusion equation, as referenced in \citet{nakshatrala2009non}, to address the transport problem using PINNs. The PDE models the transport process, providing the concentration at each spatial location $\mathbf{x}$. Let us consider a two-dimensional domain $\Omega$ with piecewise smooth boundary $\partial \Omega$. The boundary is divided into two disjoint parts: $\Gamma^\mathrm{D}$ and $\Gamma^\mathrm{N}$, such that $\Gamma^\mathrm{D} \cup \Gamma^\mathrm{N} = \partial \Omega$ and $\Gamma^\mathrm{D} \cap \Gamma^\mathrm{N} = \emptyset$ ensuring mathematical well-posedness. Here, the superscript N denotes the Neumann boundary condition, while D represents the Dirichlet boundary condition. The governing equation for second-order elliptic anisotropic diffusion is presented as \citet{nakshatrala2009non}:
%-------------------------------------------;
%  Equation: Transport governing equations  ;
%-------------------------------------------;
\begin{align}
    \label{Eqn:PINNS_Fast_Diffusion_BOM}
    -&\mathrm{div}[\mathbf{D}(\textbf{x}) \, \mathrm{grad}[c(\textbf{x})]] = f(\textbf{x})
    &&\quad \mathrm{in} \; \Omega \\
    \label{Eqn:PINNS_Fast_Diffusion_bc}
    & c(\textbf{x}) = c^{\mathrm{p}} (\textbf{x})
    &&\quad \mathrm{on} \; \Gamma^\mathrm{D} \\
    \label{Eqn:PINNS_Fast_Diffusion_bc_flux}
    -&\widehat{\mathbf{n}}(\mathbf{x}) \bullet \mathbf{D}(\textbf{x}) \, \mathrm{grad}[c(\textbf{x})] = v^{\mathrm{p}} (\textbf{x})
    &&\quad \mathrm{on} \; \Gamma^\mathrm{N}
\end{align}
where $\mathrm{D}(\textbf{x})$ is the mechanical dispersion tensor and $c(\textbf{x})$ is the concentration at spatial location $\textbf{x}$. $f(\textbf{x})$ is the volumetric source. $c^{\mathrm{p}} (\textbf{x})$ is the prescribed concentration and 
$v^{\mathrm{p}} (\textbf{x})$ is the prescribed flux at the boundary. $\widehat{\mathbf{n}}(\mathbf{x})$ is the unit outward normal to the boundary. Advection is not explicitly included in the solute transport equation, as the influence of the flow field is incorporated via the velocity-dependent mechanical dispersion tensor $\mathbf{D}(\textbf{x})$, which governs anisotropic transport behavior.

\subsection{Diffusion-reaction}
\label{SubSec:S2_DiffRxn}
\begin{figure}[hbt!]
    \centering
    \includegraphics[scale=0.85]{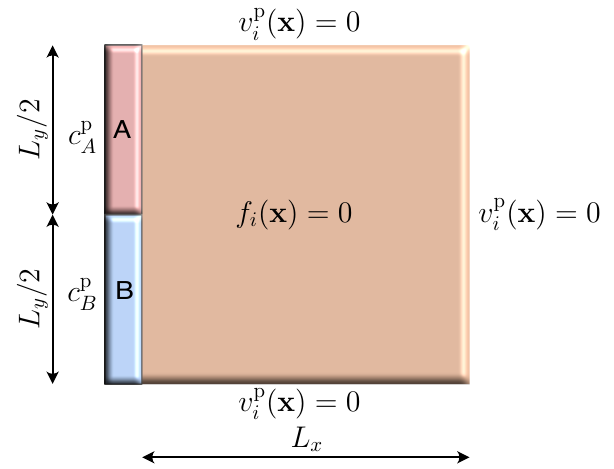}
    \caption{\textsf{Problem setup:} The figure depicts the setup of the reaction subproblem within a square domain (or reaction tank) of dimensions $L_x$ and $L_y$, and the boundary $\partial\Omega$. The left boundary is prescribed with the concentration of chemical A in the upper half and chemical B in the lower half, as illustrated, while the remaining boundaries are subjected to zero flux $v_i ^\mathrm{p}(\mathbf{x})$. $f_i(\mathbf{x})$ is the volumetric source. The subscript $i$ denotes different chemical species considered in this study: A, B, and C.}
    \label{fig:problem_setup}
\end{figure}

For the bimolecular reaction, we consider a square domain $\Omega$ (also referred to as a reaction tank) with dimensions $L_x$ and $L_y$. The boundary is denoted as $\partial \Omega$, as illustrated in Fig.~\ref{fig:problem_setup}. $\mathbf{x}$ is a spatial point in the domain, while the boundary $\partial \Omega$ is discretized into $\Gamma_i^\mathrm{D}$, where the chemical concentration is prescribed, and $\Gamma_i^\mathrm{N}$, where the chemical flux is prescribed. The subscript $i$ represents different chemical species. For the mathematical well-posedness, these satisfy the conditions $\Gamma_i^\mathrm{D} \cup \Gamma_i^\mathrm{N} = \partial \Omega$ and $\Gamma_i^\mathrm{D} \cap \Gamma_i^\mathrm{N} = \emptyset$.
%---------------------------------;
% Figure 4: Reaction Problem setup;
%---------------------------------;
For the fast bimolecular reaction, we consider two reactants, chemical A and chemical B, with concentrations denoted as $c_A$ and $c_B$, respectively. These reactants undergo a fast bimolecular reaction, resulting in the formation of chemical C, whose concentration is represented as $c_C$. Let $n_A$, $n_B$, and $n_C$ be the positive stoichiometric coefficients, such that:
\begin{align}
    n_A \, A + n_B \, B \rightarrow n_C \, C
\end{align}
The governing equation for the steady state diffusive-reactive system in a rigid porous medium with the $i-\mathrm{th}$ chemical species takes the following form:
\begin{align}
    -&\mathrm{div} \, [\mathbf{D}(\textbf{x}) \, \mathrm{grad}(c_i)] = f_i(\mathbf{x})-n_i \, r
    &&\quad \mathrm{in} \; \Omega \\
    & c_i (\mathbf{x}) = c_i ^\mathrm{p}(\mathbf{x})
    &&\quad \mathrm{on} \; \Gamma^{\mathrm{D}} \\
    -&\mathbf{D}(\textbf{x}) \; \mathrm{grad} [c_i] \bullet \widehat{\mathbf{n}}(\mathbf{x})= v_i ^\mathrm{p}(\mathbf{x})
    &&\quad \mathrm{on} \; \Gamma^{\mathrm{N}}
\end{align}
where $i = A, B, C$ and $r$ is the rate of reaction. $f_i(\mathbf{x})$ is the volumetric source, $c_i^\mathrm{P}(\mathbf{x})$ is the prescribed concentration, and $v_i ^\mathrm{p}(\mathbf{x})$ is the prescribed mass concentration flux on the boundary. $\mathbf{D}(\textbf{x})$ is the mechanical dispersion tensor calculated utilizing the velocity $\mathbf{v}$ from the flow problem such that:
\begin{align}
    \label{Fast_PINNs_diffusivity_tensor}
    \mathbf{D}(\mathbf{x}) = \alpha_\mathrm{T} \, \Vert \mathrm{\mathbf{v}} \Vert \, \mathbf{I} + \frac{(\alpha_\mathrm{L}-\alpha_\mathrm{T})}{\Vert \mathrm{\mathbf{v}} \Vert} \, \mathrm{\mathbf{v}} \otimes \mathrm{\mathbf{v}}
\end{align}
where $\mathbf{I}$ is the second-order identity tensor, and $\alpha_\mathrm{L}$ and $\alpha_\mathrm{T}$ are longitudinal and transverse dispersion, respectively \citep{pinder2006subsurface, nakshatrala2013numerical}. For simplicity, the volumetric source $f_i(\mathbf{x})$ is considered zero. On expanding for each chemicals, we get:
\begin{align}
\label{expanded}
    -&\mathrm{div} \, [\mathbf{D}(\textbf{x}) \, \mathrm{grad}(c_A)] = -n_A \, r
    &&\quad \mathrm{in} \; \Omega \\
    -&\mathrm{div} \, [\mathbf{D}(\textbf{x}) \, \mathrm{grad}(c_B)] = -n_B \, r
    &&\quad \mathrm{in} \; \Omega \\
    -&\mathrm{div} \, [\mathbf{D}(\textbf{x}) \, \mathrm{grad}(c_C)] = n_C \, r
    &&\quad \mathrm{in} \; \Omega 
\end{align}
The three coupled tensorial diffusion-reaction equations are reformulated in terms of invariants that remain unaffected by the reaction. These two invariants are expressed as follows:
\begin{align}
    \label{Invariant_1}
    & \Psi_A = c_A + \frac{n_A}{n_C} \, c_C \\
    \label{Invariant_2}
    & \Psi_B = c_B + \frac{n_B}{n_C} \, c_C
\end{align}
Now, the Eq.~\eqref{expanded} in terms of invariant $\Psi_A$ takes the following form:
\begin{align}
    \label{Fast_PINNs_invariant1}
    -&\mathrm{div} \, [\mathbf{D}(\textbf{x}) \, \mathrm{grad}(\Psi_A)] = 0
    &&\quad \mathrm{in} \; \Omega \\
    & \Psi_A^{\mathrm{p}} = c_A^{\mathrm{p}} + \frac{n_A}{n_C} \, c_C^{\mathrm{p}} 
    &&\quad \mathrm{on} \; \Gamma^{\mathrm{D}}
\end{align}
Likewise,
\begin{align}
    \label{Fast_PINNs_invariant2}
    -&\mathrm{div} \, [\mathbf{D}(\textbf{x}) \, \mathrm{grad}(\Psi_B)] = 0
    &&\quad \mathrm{in} \; \Omega\\
    \label{Fast_PINNs_invariant2.}
    & \Psi_B^{\mathrm{p}} = c_B^{\mathrm{p}} + \frac{n_B}{n_C} \, c_C^{\mathrm{p}} 
    &&\quad \mathrm{on} \; \Gamma^{\mathrm{D}}
\end{align}
The above equations Eqs.~\eqref{Fast_PINNs_invariant1}-\eqref{Fast_PINNs_invariant2.} represent an initial boundary value problem where the invariants $\Psi_A$ and $\Psi_B$ can be solved throughout the domain.

%============================================;
%  Subsubsection: Fast bimolecular reaction  ;
%============================================;
\subsubsection{Fast bimolecular reaction} 
A fast reaction refers to a scenario where the reaction timescale is significantly shorter than the diffusion timescales of the chemical species. Consequently, for fast bimolecular reactions, it is a reasonable approximation to assume that species A and B cannot coexist at the same location $\mathbf{x}$ at any given time. This assumption becomes exact in the limit of infinitely fast or instantaneous reactions. To determine the concentration of each chemical in such fast reactions, we rely on expressions for invariants. On multiplying Eq.~\eqref{Invariant_1} with $n_B$ and ~\eqref{Invariant_2} with $n_A$ and subtracting, we get:
\begin{align}
    \label{ca_cb}
    n_B \, \Psi_A - n_A \, \Psi_B = n_B \, c_A - n_A \, c_B
\end{align}
Since we are considering a fast reaction, if the expression in Eq.~\eqref{ca_cb} is greater than zero, it implies that there exists a concentration $c_A$ for which $c_B$ cannot coexist, i.e.,
\begin{align}
    & c_A = \frac{n_B \, \Psi_A - n_A \, \Psi_B}{n_B} \\
    & c_B = 0
\end{align}
Likewise, if the expression in Eq.~\eqref{ca_cb} is less than zero, it implies there exists a concentration $c_B$ where $c_A$ cannot coexist, i.e.,
\begin{align}
    & c_B = - \frac{n_B \, \Psi_A - n_A \, \Psi_B}{n_A} \\
    & c_A = 0
\end{align}
Therefore, we utilize maximum and minimum principles to evaluate the concentrations $c_A$ and $c_B$ as following:
\begin{align}
    \label{Fast_PINNs_cA}
    c_A = \frac{1}{n_B} \, \mathrm{max} \big[n_B \, \Psi_A - n_A \, \Psi_B, \, 0\big] \\
    \label{Fast_PINNs_cB}
    c_B = \frac{1}{n_A} \, \mathrm{min} \big[n_B \, \Psi_A - n_A \, \Psi_B, \, 0\big] 
\end{align}
Finally, $c_C$ is obtained using Eq.~\eqref{Invariant_1}:
\begin{align}
    \label{Fast_PINNs_cC}
    c_C = \frac{n_C}{n_A} \, (\Psi_A - c_A)
\end{align}
After obtaining the invariants from the PINNs training, additional mathematical computations are performed using Eqs.~\eqref{Fast_PINNs_cA}--\eqref{Fast_PINNs_cC} to determine the concentrations of chemicals A, B and C, respectively.

%============================================================;
%  Subsubsection: Applicability to critical mineral science  ;
%============================================================;
\subsubsection{Applicability of fast bimolecular reactions to critical mineral science} 
When reaction rates greatly exceed advective or diffusive transport rates, local chemical equilibrium can be assumed throughout the flow field. However, when kinetic reactions and transport occur on comparable timescales, sharp reaction fronts---such as mineral dissolution waves---may develop. For processes like adsorption, precipitation, and redox reactions, kinetic limitations can strongly influence the extraction of rare earth elements (REEs). In such cases, it is necessary to explicitly model both fast reactions and finite-rate kinetics to accurately capture plume evolution and mineralogical zoning \citep{steefel2009fluid, zhou2023geochemical}.

Many critical minerals (e.g., REEs, Co, Ni) are recovered through acid- or ligand-mediated \emph{in situ} leaching. The dissolution of mineral grains by complexing agents often follows fast bimolecular kinetics, such as:
\begin{align}
    \mathrm{Mineral}_{(s)} \;+\; \mathrm{H}^+_{(\mathrm{aq})} \;\longrightarrow\; \mathrm{Metal}^{n+}_{(\mathrm{aq})} \;+\; \cdots \\
    \mathrm{Metal(CO_3)_2}_{(s)} \;+\; 2\,\mathrm{H}^{+}_{(\mathrm{aq})}
    \;\longrightarrow\;
    \mathrm{Metal}^{2+}_{(\mathrm{aq})} \;+\; 2\,\mathrm{HCO}_{3 \, (\mathrm{aq})}^{-}
\end{align}
When $\mathrm{H}^+$ ions collide with reactive sites on the mineral surface, they can break bonds within the solid structure, releasing metal ions into solution. For many oxides and silicates, this process is surface-controlled and follows a rate law of the form
\begin{align}
    r = k_{\mathrm{surf}} \, a^{n}_{\mathrm{H}^{+}} 
\end{align}
where $k_{\mathrm{surf}}$ is the surface-reaction rate constant and $a_{\mathrm{H}^{+}}$ is the proton activity. This rate law is valid primarily under low pH conditions \citep{brantley2008analysis} because $n$ is often unity for simple proton-mediated dissolution. In such scenarios, the reaction is formally bimolecular (involving one $\mathrm{H}^{+}$ and one reactive surface site). The rate can be very high, reaching near-equilibrium on laboratory timescales for minerals such as calcite or barite. When $k_{\mathrm{surf}}$ is large relative to advective or diffusive transport rates, local equilibrium can be assumed, allowing the use of algebraic constraints in reactive-transport models. However, when reaction and transport rates are comparable, reaction fronts may emerge, with their shape and velocity governed by the interplay between kinetics and transport processes.

Recognizing that certain reactions operate in the fast limit allows the use of algebraic equilibrium constraints (e.g., \eqref{Invariant_1}–\eqref{Invariant_2}) in place of stiff ordinary differential equations (ODEs).
This simplification reduces model complexity, making it easier for PINNs to learn underlying dynamics.
In critical mineral recovery workflows, assuming fast kinetics enables PINNs to focus on capturing the evolution of chemical plumes in heterogeneous flow fields. By avoiding the need to resolve stiff reaction terms, these models require fewer computational resources, converge more quickly, and yield more interpretable results. In practical settings---such as heap leach pads or in situ recovery wells---real-time monitoring data (e.g., fluid chemistry) can be integrated into reactive transport models. In this context, fast bimolecular kinetics modeled via PINNs enable real-time updates of plume dynamics and recovery predictions. Accurately resolving reactive mixing and front propagation improves the ability of PINNs to capture sharp reaction fronts and localize geochemical changes in heterogeneous media. Under fast bimolecular kinetics, such behavior can be leveraged in reactive transport simulations to reduce computational stiffness and improve predictive fidelity in systems where reaction and transport timescales are comparable. 

A summary of workflow steps is as follows:
\begin{enumerate}
  \item \textbf{Input \& setup:} Specify domain/boundary conditions, heterogeneous permeability $K(\mathbf{x})$, porosity $\phi$, fluid properties, stoichiometry $(n_A,n_B,n_C)$, and dispersion parameters $(\alpha_L,\alpha_T,D_m)$.

  \item \textbf{Flow PINN (Darcy/mixed form):} Train a PINN to solve for pressure $p(\mathbf{x})$ and Darcy velocity $\mathbf{v}(\mathbf{x})$ using a PDE-residual loss plus boundary-condition losses; validate with patch tests. The trained velocity field $\mathbf{v}(\mathbf{x})$ is then fixed and passed to the transport step.

  \item \textbf{Compute anisotropic dispersion tensor:} Construct
  \[
    \mathbf{D}(\mathbf{x}) \;=\; D_m \mathbf{I} \;+\; \mathbf{D}_{\mathrm{disp}}\!\big(\mathbf{v}(\mathbf{x}),\alpha_L,\alpha_T\big),
  \]
  and, if used, the effective tensor $\mathbf{D}(\mathbf{x})/\phi$.

  \item \textbf{Transport--reaction via invariants (reaction part):} Instead of training a network directly on stiff reaction source terms, solve two \textbf{reaction-invariant} transport equations for $\Psi_A(\mathbf{x},t)$ and $\Psi_B(\mathbf{x},t)$. These invariants are constructed so that the fast bimolecular reaction cancels from the governing equations, leaving uncoupled diffusion equations (implemented here through anisotropic diffusion with $\mathbf{D}(\mathbf{x})$).

  \item \textbf{Algebraic fast-reaction closure (reconstruction):} Recover species concentrations from invariants using stoichiometry and non-negativity/mutual-exclusion rules, e.g.,
  \[
    c_A=\max\!\Big(\frac{\Psi_A}{n_A}-\frac{\Psi_B}{n_B},0\Big),\qquad
    c_B=\max\!\Big(\frac{\Psi_B}{n_B}-\frac{\Psi_A}{n_A},0\Big),\qquad
    c_C=\min\!\Big(\frac{\Psi_A}{n_A},\frac{\Psi_B}{n_B}\Big),
  \]
  (or equivalent piecewise forms). This step concerns the ``ML implementation for the reaction'' in our approach: the network learns the invariant transport fields, while the reaction itself is enforced via an explicit, physics-based algebraic projection that is consistent with the instantaneous-reaction limit.

  \item \textbf{Outputs \& checks:} Generate $c_A,c_B,c_C$ fields, verify non-negativity, and compare plume metrics (e.g., front location and product evolution) against reference solutions.
\end{enumerate}

\subsubsection*{Scope and extensibility} 
In this proof-of-concept, we restrict attention to a single irreversible bimolecular reaction $A + B \rightarrow C$ to focus on issues specific to fast reactions. However, we note that the PINN formulation is not limited to a single reaction: multi-reaction networks can be handled by (i) augmenting the state vector to include all relevant species, (ii) adding the corresponding mass-action source terms to the PDE residuals, and (iii) enforcing algebraic equilibrium constraints as additional residuals (or via Lagrange multipliers) when present. Stoichiometric conservation and total-mass constraints are promoted through integral or boundary-condition loss terms. A systematic evaluation of mixed kinetic/equilibrium systems and larger reaction networks is an important next step and is planned for future work.

%*********************************************;
%                                             ;
%  NAME                                       ;
%    S3_PINNS_Fast_Flow.tex                   ;
%                                             ;
%*********************************************;
\section{VERIFICATION OF FLOW SUBPROBLEM}
\label{Sec:S3_PINNS_Fast_flow}
\begin{figure}[hbt!]
    \centering
    \includegraphics[scale = 0.58]{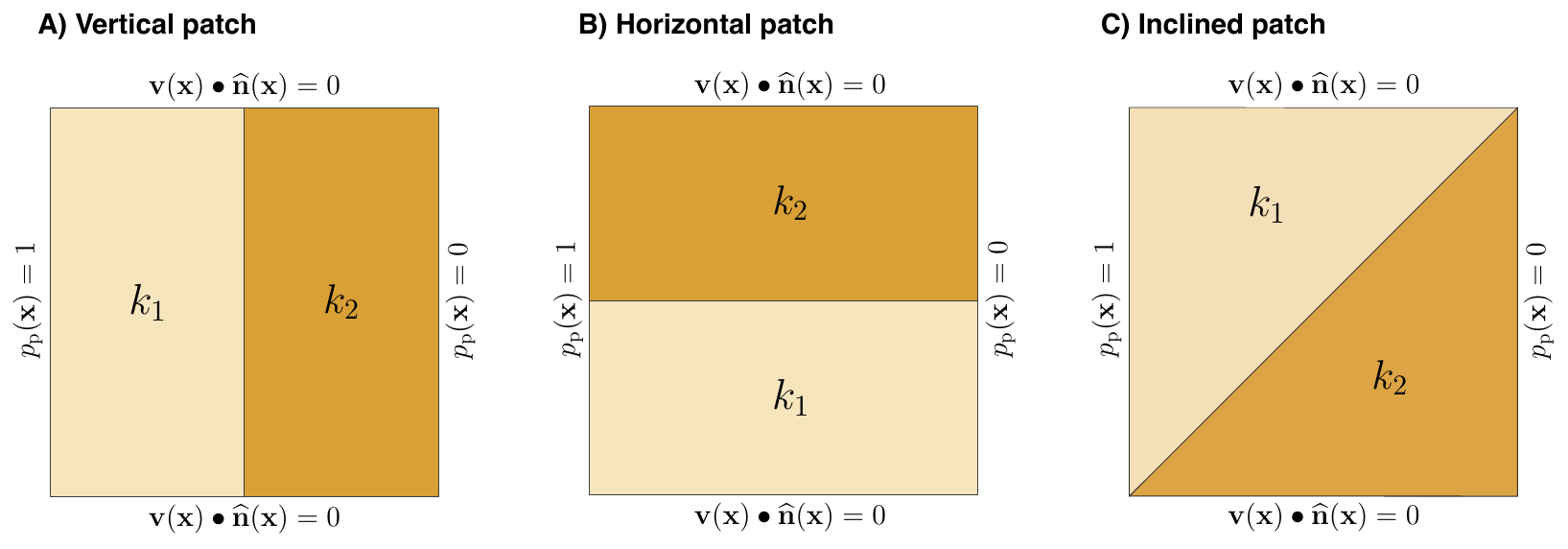}
    \caption{\textsf{Patch test:} The figure illustrates the three patch tests: vertical, horizontal, and inclined. In all cases, pressure is prescribed at the left and right boundaries, while the top and bottom boundaries have zero velocity. The permeability coefficients, $k_1$ and $k_2$ are assigned according to the specific configuration of each test setup.}
    \label{fig:PINNs_fast_patch_setup}
\end{figure}

Before proceeding to the fast bimolecular reaction, it is essential to test the developed PINNs code for the flow problem. Given that the study focuses on subsurface applications, where heterogeneous material properties are commonly encountered, evaluating the code's performance in a heterogeneous setting is a crucial step in assessing its applicability to both homogeneous and heterogeneous cases. To achieve this, we introduce three patch tests: vertical, horizontal, and inclined, as illustrated in Fig.~\ref{fig:PINNs_fast_patch_setup}.

In all cases, a unit pressure is prescribed at the left boundary and zero pressure at the right boundary, while the top and bottom boundaries have zero velocity. The domain consists of two permeability values, $k_1$ and $k_2$, which are arranged differently depending on the patch test. In the vertical patch test, the left half of the domain has a permeability of $k_1 = 1$, while the right half has $k_2 = 10$. In the horizontal patch test, the bottom half of the domain has 
$k_1 = 1$, and the top half has $k_2 = 10$. Similarly, in the inclined patch test, the top diagonal region has a permeability of 
$k_1 = 1$, while the bottom diagonal region has $k_2 = 10$.

For all three cases, the boundary value equation presented in Eq.~\eqref{Eqn:PINNS_Fast_BoLM} is used to train the PINNs model and solve for pressure and velocity at each spatial point, a total of $100\times100$ for each case. The solution is then tested and validated against the finite element method (FEM) solution obtained utilizing fluid flow models---second-order Lagrange elements on a triangular mesh with a standard continuous Galerkin formulation---available in \citet{COMSOL}.

\subsubsection*{Assumption on permeability} 
To reduce confounding factors and establish a baseline, we use constant (or piecewise constant) permeability fields in all tests. In realistic settings, permeability varies spatially and can evolve in time due to mineralization/dissolution and porosity changes. Within the same PINNs framework, heterogeneous or time-dependent permeability $\mathbf{K}(\mathbf{x},t)$ can be: (i) provided as an external field and incorporated directly into the Darcy residuals; or (ii) inferred from data by parameterizing $\mathbf{K}(\mathbf{x},t)$ (e.g., a deep neural network or low-dimensional basis) and training jointly with state variables using pressure/flux observations alongside Darcy's law and continuity in the loss. Constitutive links such as $\mathbf{K}=\mathbf{K}(\phi(c_i))$ can be included straightforwardly. $\phi$ is the porosity. Evaluating performance under heterogeneous and evolving $\mathbf{K}$ is a natural extension that we will address in future work.

%------------------------;
%  Figure 4: Patch test  ;
%------------------------;

\subsection{Vertical patch}
\begin{figure}[hbt!]
    \centering
    \includegraphics[scale = 0.68]{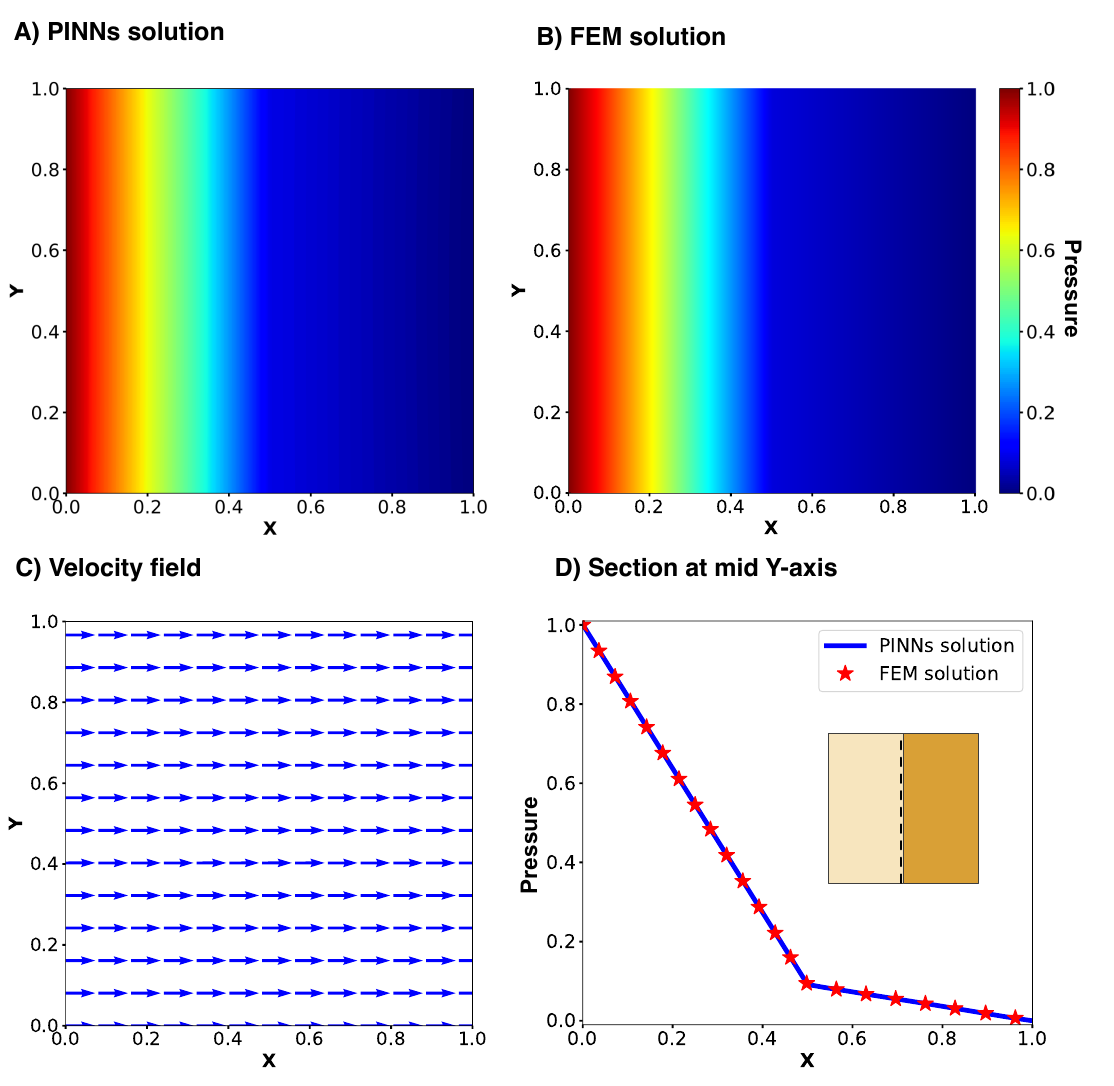}
    \caption{\textsf{Vertical patch:} The figure presents the results of the vertical patch test. The pressure profiles obtained from PINNs and FEM are comparable and appear similar in A) and B). The velocity field shown in C) depicts a uniform flow from left to right. Finally, D) displays the pressure distribution along the length at the mid y-axis, confirming that the PINNs solution aligns closely with the numerical solution.}
    \label{fig:PINNs_fast_vertical_patch}
\end{figure}

First, PINNs is applied to the vertical patch test, where the permeability coefficient is $k_1 = 1$ on the left half and $k_1 = 10$ on the right half. The heterogeneous flow problem is solved using PINNs and validated against the numerical solution, as shown in Fig.~\ref{fig:PINNs_fast_vertical_patch}. The pressure profile is first examined, revealing a strong agreement between the PINNs and numerical solutions. The maximum pressure occurs at the prescribed left boundary, while the minimum pressure is zero at the right boundary, satisfying the boundary conditions. To further analyze the results, the velocity field is also plotted, showing a uniform flow from left to right. For a clearer comparison, a section along the mid y-axis is plotted, highlighting the agreement between the two solutions. In both cases, a sharp change in slope is observed at the interface where permeability changes. At each point, the PINNs and numerical solutions are close, validating the PINNs approach against the numerical solution.
%----------------------------;
%  Figure 5: Vertical patch  ;
%----------------------------;
\subsection{Horizontal patch}
\begin{figure}[hbt!]
    \centering
    \includegraphics[scale = 0.68]{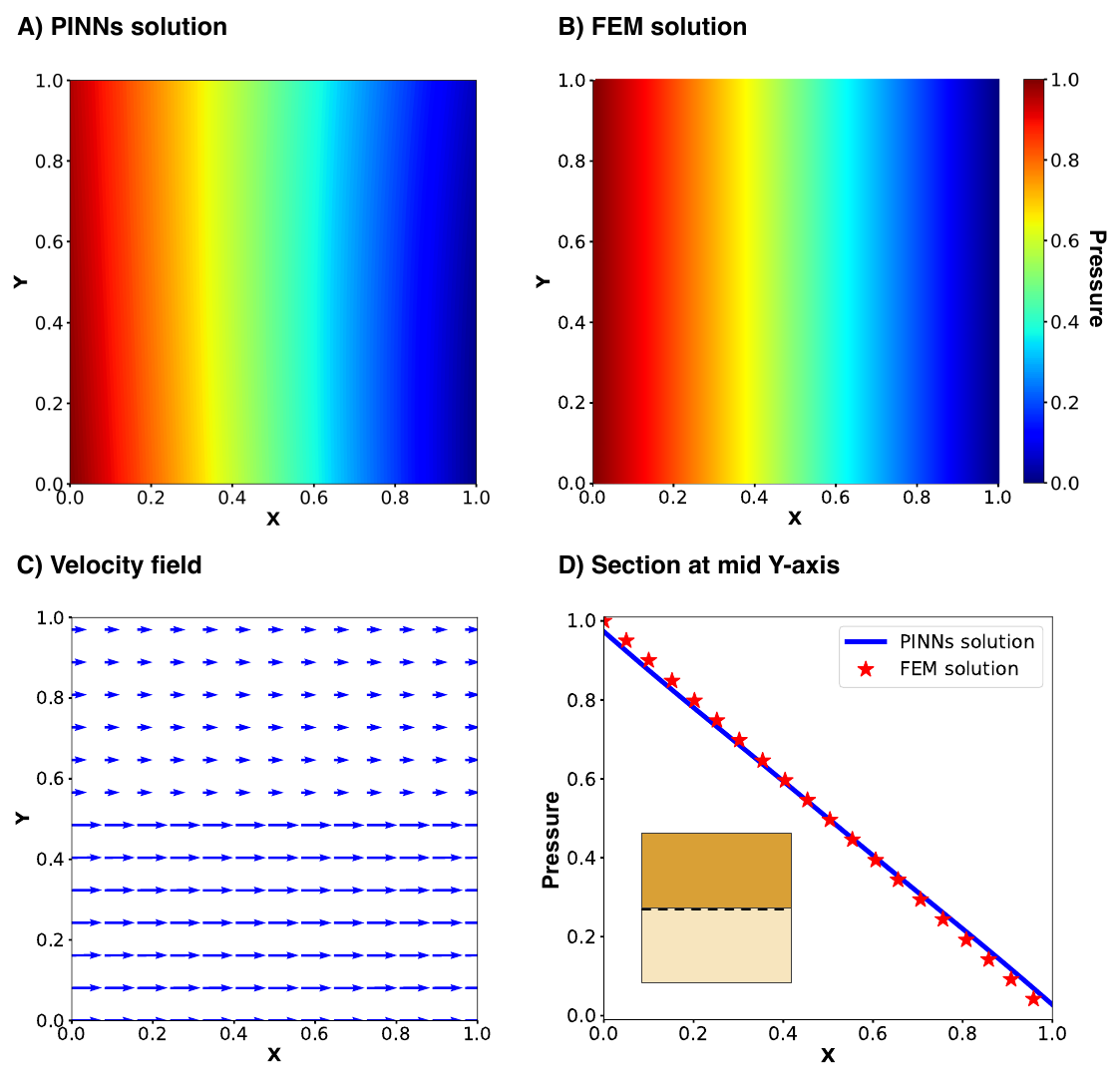}
    \caption{\textsf{Horizontal patch:} The figure presents the results of the horizontal patch test. The pressure profiles obtained from PINNs and FEM are shown in A) and B). The velocity field, which varies in magnitude in different regions due to the heterogeneous permeability distribution, is depicted in C). In D), the pressure distribution along the x-axis at the mid y-axis is shown for both the PINNs and numerical solutions. The results from both methods are comparable, showing an almost identical pressure distribution.}
    \label{fig:PINNs_fast_horizontal_patch}
\end{figure}

Second, a horizontal patch test is conducted using the developed PINNs code. In this test, the permeability is $k_1 = 1$ in the bottom half of the domain and $k_2 = 10$
in the top half. The results of the horizontal patch test are shown in Fig. ~\ref{fig:PINNs_fast_horizontal_patch}. The pressure profiles obtained from PINNs and FEM are comparable and closely match, satisfying the prescribed boundary conditions. The velocity field indicates a uniform flow from left to right, however the flow magnitude differs between top and bottom regions as reflected by the arrow sizes. Also, the pressure distribution along the x-axis at the mid y-axis is plotted, demonstrating a good agreement between the PINNs and numerical solutions, with only minor to negligible variations.
%------------------------------;
%  Figure 6: Horizontal patch  ;
%------------------------------;
\subsection{Inclined patch}
\begin{figure}[hbt!]
    \centering
    \includegraphics[scale = 0.68]{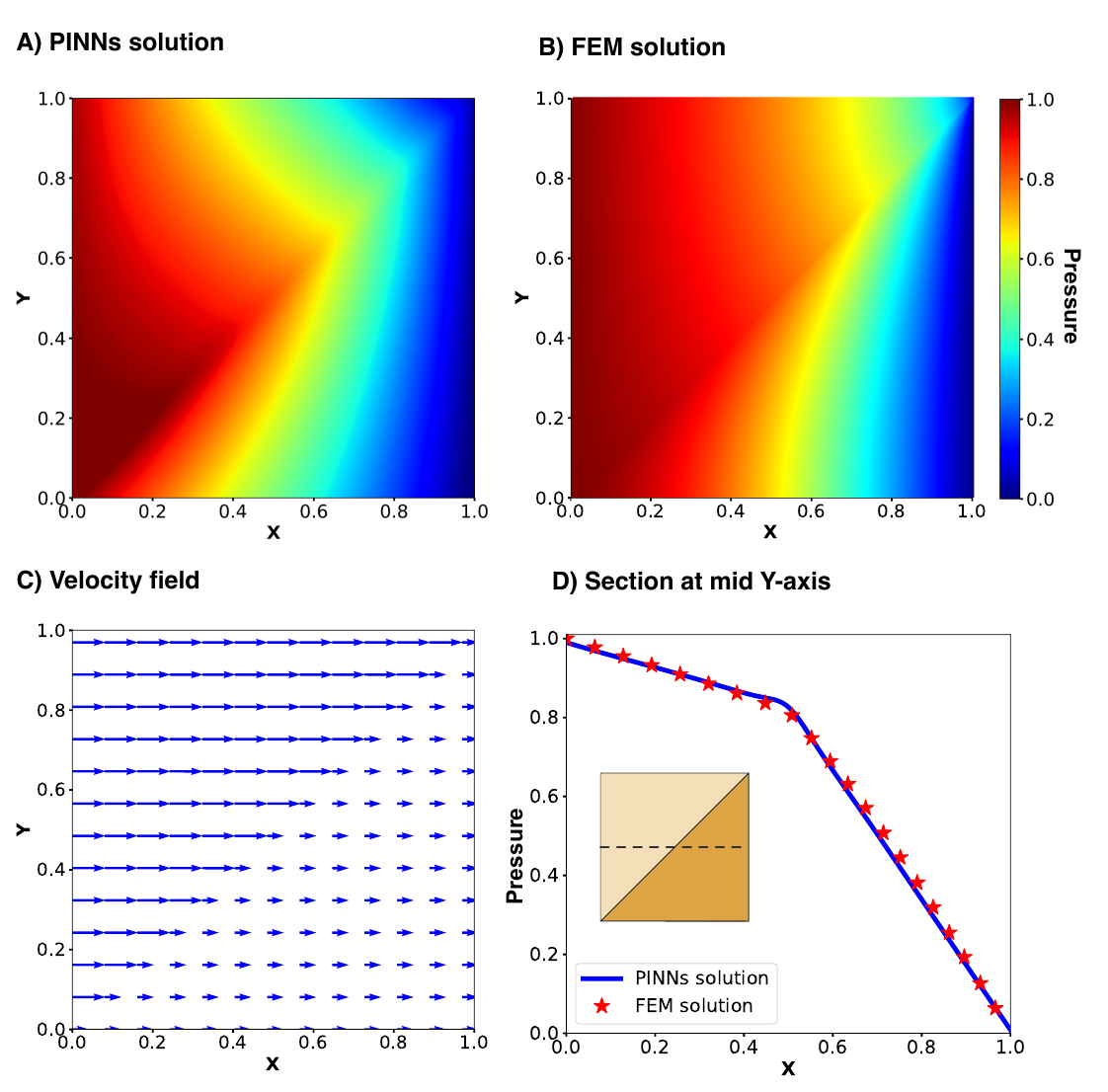}
    \caption{\textsf{Inclined patch:} The figure illustrates the inclined patch test, showing comparable pressure profiles from PINNs and FEM in A) and B), respectively. The velocity field in C) depicts a uniform flow, with variations in magnitude due to permeability heterogeneity, as indicated by the arrow sizes. Finally, the pressure distribution along the x-axis at the mid y-axis, shown in D), demonstrates a strong agreement between the PINNs and numerical solutions.}
    \label{fig:PINNs_fast_inclined_patch}
\end{figure}

Third, an inclined patch test is conducted using the PINNs code and validated against the FEM solution, as shown in Fig. ~\ref{fig:PINNs_fast_inclined_patch}. In this test, the top diagonal half of the domain has a permeability of $k_1 = 1$, while the bottom diagonal half has $k_2 = 10$. The pressure profiles from both PINNs and FEM capture the impact of heterogeneous permeability at the interface while satisfying the boundary conditions. The velocity field indicates a uniform flow, though its magnitude varies due to the permeability differences in the diagonal regions. Finally, the pressure distribution along the x-axis for both PINNs and numerical solutions is close enough, confirming the accuracy and validity of the obtained results.
%----------------------------;
%  Figure 7: Inclined patch  ;
%----------------------------;

All the results combined together in Figs. ~\ref{fig:PINNs_fast_vertical_patch}--\ref{fig:PINNs_fast_inclined_patch} complete the three patch tests, each highlighting different scenarios of permeability heterogeneity. The error analysis and convergence plots are presented in the supplementary material to provide robustness in the solution obtained. The developed PINNs code successfully captures the heterogeneity in all cases, demonstrating its effectiveness for such flow problems. The results closely match those obtained using FEM, confirming the accuracy of the approach. Thus, this section provides a validation of the PINNs' capability to solve these types of flow problems, hence answering the first scientific question. 

%*********************************************;
%                                             ;
%  NAME                                       ;
%    S4_PINNS_Fast_Transport.tex              ;
%                                             ;
%*********************************************;
\section{NON-NEGATIVE SOLUTIONS UNDER THE TRANSPORT PROBLEMS}
\label{Sec:S4_PINNS_Fast_Transport}

\subsection{Maximum principle} A key property of diffusion-type equations is that their solutions remain non-negative and adhere to the maximum principle. Physically, the concentration in a diffusion problem is non-negative. If $c(\mathbf{x})$ represents the concentration over the domain i.e., $c(\mathbf{x}) \in C^2(\Omega) \cap C(\overline{\Omega})$ where $C^2(\Omega)$ denotes the set of twice continuously differentiable functions, and satisfies the diffusion equation in Eq. ~\eqref{Eqn:PINNS_Fast_Diffusion_BOM}, then the minimum concentration occurs at the boundary, given by:
\begin{align}
    \label{Fast_PINNs_maximum_principle}
    \min_{\mathbf{x}\in\overline{\Omega}} \; c(\mathbf{x}) = \min_{\mathbf{x}\in{\partial \Omega}} c(\mathbf{x})
\end{align}
However, in previous studies, various numerical methods, such as the finite element method \citep{braess2001finite, barrenechea2024finite}, finite difference method \citep{morton2005numerical, strikwerda2004finite}, and finite volume method \citep{droniou2006mixed}, often struggle or fail to strictly satisfy the maximum principle when solving transport problems with high accuracy. Studies have explored enforcing the maximum principle \citep{evans2009enforcement,herrera2006positive}, optimization techniques \citep{nakshatrala2009non} and have also employed nonlinear discretization  to improve adherence to the maximum principle while achieving reasonably accurate solutions \citep{barrenechea2024finite}. Given the challenges that popular numerical methods face in maintaining non-negative solutions for diffusion-type problems, it is worthwhile to examine the performance of PINNs and determine whether this novel approach can uphold the maximum principle or not.
\subsection{Problem setup} We consider a square domain with a unit length $a = 1$ featuring a centrally located square hole with a side length of 0.2 units. The inner boundary denoted by $\Gamma_\mathrm{in}$ is prescribed with a maximum concentration $c^{\mathrm{p}} (\textbf{x}) = 1 \in \Gamma_\mathrm{in}$ while the outer boundary denoted by $\Gamma_\mathrm{out}$ is set to a minimum concentration of $c^{\mathrm{p}} (\textbf{x}) = 0 \in \Gamma_\mathrm{out}$. The forcing function $f(\mathbf{x})$ is taken as zero. Given that the study solves the transport problem in a 2D anisotropic medium, the mechanical dispersion tensor $\mathrm{D}(\textbf{x})$ is defined as:
\[
   \mathrm{D}=
  \left[ {\begin{array}{cc}
   \cos({\theta}) & \sin({\theta}) \\
   -\sin({\theta}) & \cos({\theta}) \\
  \end{array} } \right] \; \left[ {\begin{array}{cc} 
   \lambda_1 & 0 \\
   0 & \lambda_2 \\
  \end{array} } \right] \; \left[ {\begin{array}{cc}
   \cos({\theta}) & -\sin({\theta}) \\
   \sin({\theta}) & \cos({\theta}) \\
  \end{array} } \right]
\]
where $\theta = \pi/6$, the principal components of the anisotropic diffusion $\lambda_1 = 10000$ and $\lambda_2 = 1$ for extreme anisotropy. With this configuration, the governing equations in Eqs. ~\eqref{Eqn:PINNS_Fast_Diffusion_BOM}--\eqref{Eqn:PINNS_Fast_Diffusion_bc_flux} are solved over the square plate with a hole.
\subsection{Numerical results}
A two-dimensional diffusion problem in an anisotropic medium is solved using two approaches: the finite element method (FEM) in \citet{COMSOL}---second-order Lagrange elements on a triangular mesh with a standard continuous Galerkin formulation---and the novel PINNs-based method introduced in this study. The problem is formulated on a square domain with a central hole, where the concentration is computed based on the governing equations and problem setup outlined in the previous subsections. A key aspect of the analysis is evaluating whether each approach adheres to the maximum principle by ensuring non-negative solutions.

\begin{figure}[hbt!]
    \centering
    \includegraphics[scale=0.75]{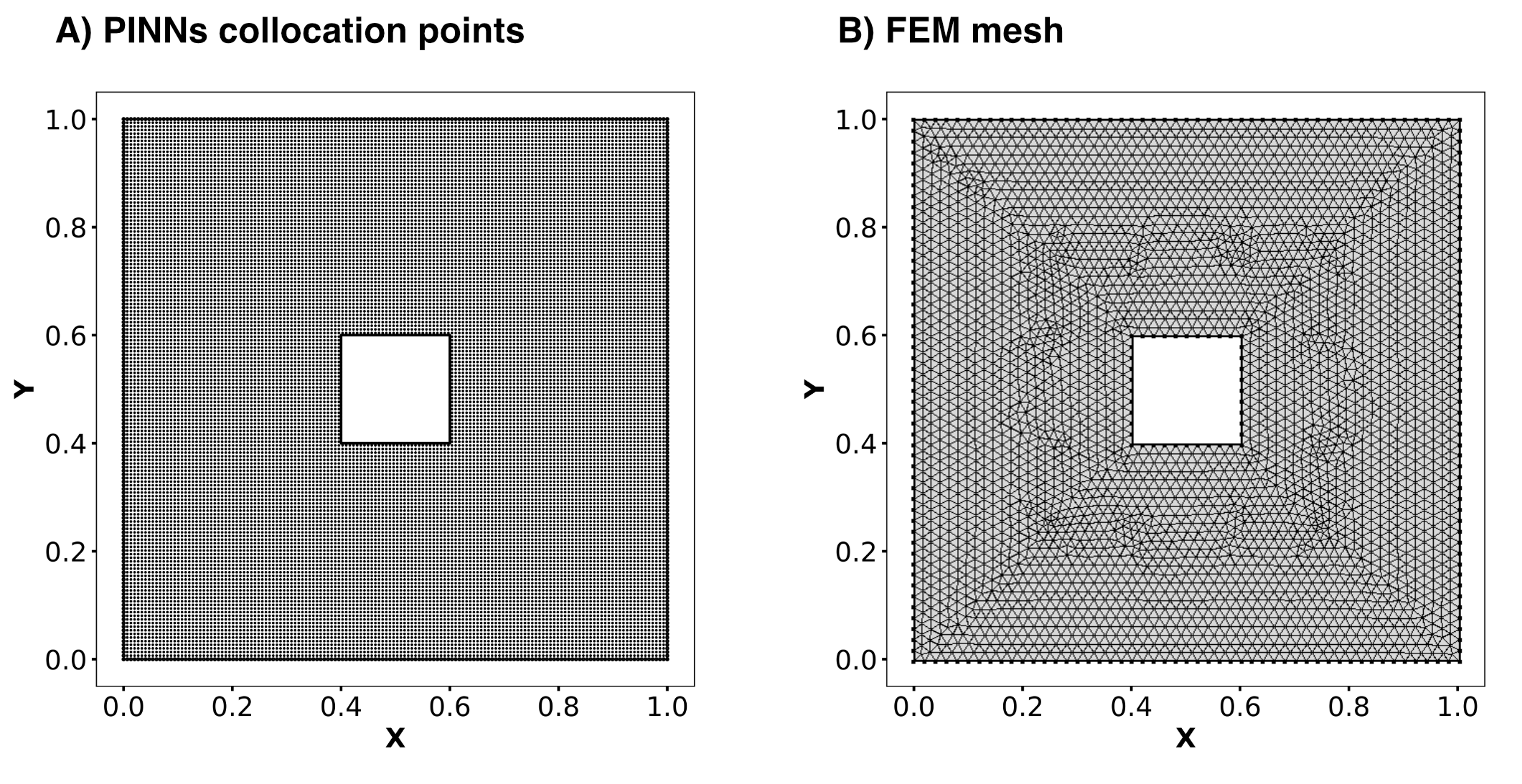}
    \caption{\textsf{Discretization strategies:} The figure illustrates two different domain discretization approaches. (A) depicts the discretization of a square plate with a hole using collocation points, a technique employed in the PINNs framework to enable a mesh-free environment. (B) illustrates well-centered triangular meshing, which is used to solve the transport problem via the finite element method.}
    \label{fig:Fast_PINNs_transport_mesh}
\end{figure}

First, the domain is discretized using two different approaches: collocation points and meshing, as illustrated in Fig. ~\ref{fig:Fast_PINNs_transport_mesh}. For solving with PINNs, discretization is performed using collocation points, with 150 points placed along each side, resulting in a total of 22,500 points within the domain. This method eliminates the need for traditional mesh generation (e.g., elements and connectivity), offering flexibility in selecting collocation points based on geometry and solution requirements. A triangular mesh with well-centered elements—where the centroid lies close to the geometric center of each triangle—is employed to improve accuracy for anisotropic diffusion problems. The mesh is refined using COMSOL’s ``extra-fine'' setting, yielding minimum element size of approximately 0.0015 unit. However, further reducing the mesh size can increase computational costs, making it less practical for large-scale simulations.
%-----------------------------------;
%  Figure 8: Domain discretization  ;
%-----------------------------------;

\begin{figure}[hbt!]
    \centering
    \includegraphics[scale = 0.7]{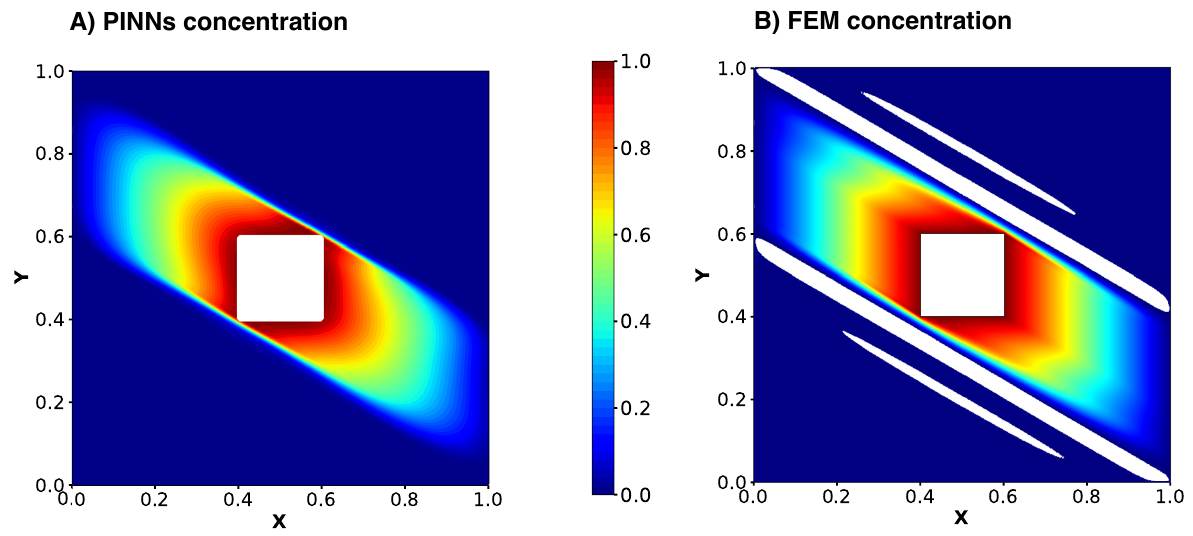}
    \caption{\textsf{Concentration profile:} The figure presents the concentration profiles obtained using PINNs and FEM in (A) and (B), respectively. Both approaches successfully solve the anisotropic diffusion problem. A key observation is that PINNs adhere to the maximum principle, yielding only non-negative concentrations across the domain, with the highest values near the inner boundaries and the lowest near the outer boundaries. In contrast, FEM results exhibit negative concentrations, represented by white regions, indicating a violation of the maximum principle.}
    \label{fig:Fast_PINNs_transport_result}
\end{figure}

The diffusion problem in Eqs.~\eqref{Eqn:PINNS_Fast_Diffusion_BOM}--\eqref{Eqn:PINNS_Fast_Diffusion_bc_flux} are solved for concentration using both PINNs and FEM, as illustrated in Fig. ~\ref{fig:Fast_PINNs_transport_result}. Both methods effectively capture anisotropic diffusion and produce comparable concentration profiles across the domain. The maximum concentration occurs at the inner boundary $c(\textbf{x}) = 1$, as prescribed . However, according to the maximum principle in Eq. ~\eqref{Fast_PINNs_maximum_principle}, the minimum concentration should be the minimum value occuring at the boundary, which is $c^{\mathrm{p}} (\textbf{x}) = 0$ prescribed at the outer boundary in the problem setup.
%--------------------------;
%  Figure 9: Transport NR  ;
%--------------------------;
The PINNs solution inherently satisfies the maximum principle, ensuring all concentration values remain non-negative, with the lowest concentration $c(\textbf{x}) = 0$. In contrast, the FEM approach fails to maintain the maximum principle, resulting in negative concentrations, which are represented by white color in Fig. ~\ref{fig:Fast_PINNs_transport_result}B). To enforce the maximum principle in FEM, additional constraints and modifications are required \citep{mudunuru2016enforcing, barrenechea2024finite}.

On the other hand, PINNs produced solutions that adhered to the maximum principle without requiring additional constraints or modifications to the loss function or model architecture, responding to the second science question. Unlike classical discretization methods, PINNs approximate the solution as a smooth global function and enforce the governing equations pointwise at collocation points via the loss function. This strong-form enforcement, together with the absence of element-wise polynomial interpolation, mitigates spurious oscillations commonly associated with Gibbs phenomena. Consequently, non-physical negative concentrations are implicitly penalized during training, allowing the PINNs solution to adhere to the maximum principle without additional constraints. This highlights the advantage of PINNs in solving diffusion problems while strictly maintaining physical principles, making them a powerful and efficient tool for such applications.

%*********************************************;
%                                             ;
%  NAME                                       ;
%    S5_PINNS_Fast_NR.tex                     ;
%                                             ;
%*********************************************;
\section{NUMERICAL RESULTS: FAST BIMOLECULAR REACTIONS}
\label{Sec:S5_PINNS_Fast_NR}
The steady-state diffusion-reaction equations in section \ref{SubSec:S2_DiffRxn} representing a fast irreversible bimolecular reaction are solved to determine the concentration of reactants and products while also investigating plume formation within the domain using PINNs. The problem setup follows the configuration shown in Fig. ~\ref{fig:problem_setup}. The prescribed concentration of species A on the top half of the left boundary is set as $c_A^\mathrm{p} = 1$, while species B is prescribed on the bottom boundary with $c_B^\mathrm{p} = 1$. When species A is prescribed, the concentrations of B and C are set to zero, and vice versa for species B. The stoichiometric coefficients are chosen as $n_A = 1$, $n_B = 2$ and $n_C = 1$. The mechanical dispersion tensor is calculated utilizing Eq. ~\eqref{Fast_PINNs_diffusivity_tensor}. To simulate highly anisotropic transport and challenge the PINNs framework, we used a strongly directional dispersion ratio where the parameters are set to $\alpha_\mathrm{L} = 1$ and $\alpha_\mathrm{T} = 10^{-5}$, with $\mathbf{I}$ representing the identity matrix. Two different velocity fields $\mathbf{v}$ are considered to examine the fast bimolecular reactions and plume formation using PINNs and evaluate its effectiveness in solving such diffusion-reaction problems.

\subsection{Uniform flow velocity}
\begin{figure}[hbt!]
    \centering
    \includegraphics[scale = 0.65]{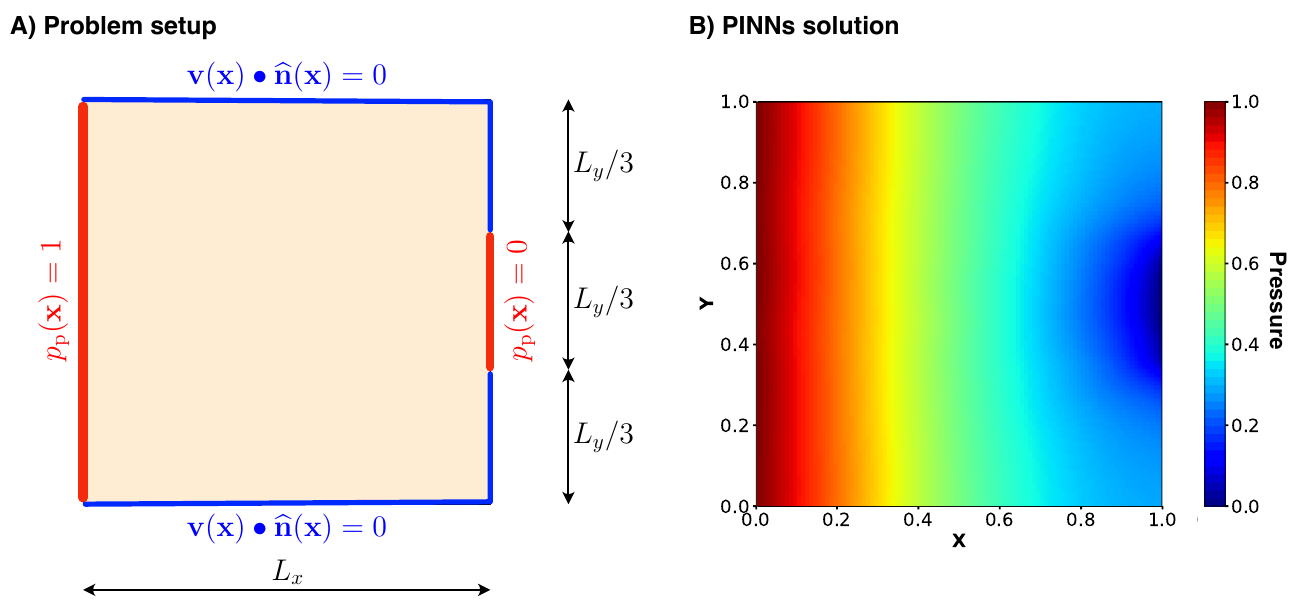}
    \caption{\textsf{Problem setup and pressure distribution:} The figure illustrates the problem setup for the flow subproblem used in the chemical reaction in A). The right boundary is assigned a unit pressure, while the middle third of the right boundary is set to zero pressure. The remaining boundaries are prescribed with zero velocity. The flow problem is solved utilizing PINNs to solve for pressure at each collocation point as shown in B). The highest pressure is observed at the left boundary, while the lowest occurs at the middle section of the right boundary, satisfying the boundary conditions defined in the problem setup. }
    \label{fig:Fast_PINNs_flow_problem_setup_pressure}
\end{figure}

\begin{figure}[hbt!]
    \centering
    \includegraphics[scale=0.65]{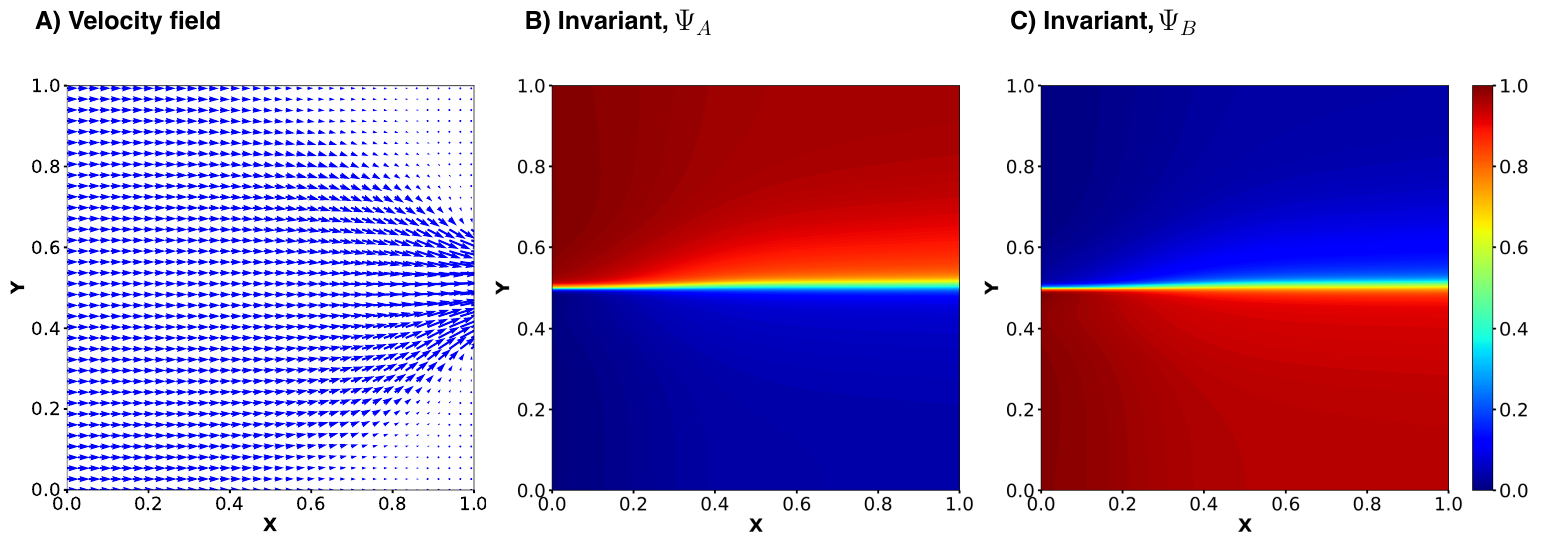}
    \caption{\textsf{Velocity field and invariants:} The figure presents the velocity field obtained from the flow subproblem in (A). The velocity field is uniform, with noticeable flow occurring only in the middle third of the right boundary, consistent with the prescribed boundary conditions. The invariants are shown in (B) and (C), where the values of $\Psi_A$ and $\Psi_B$ range from 0 to 1, as expected. This confirms that the diffusion remains bounded and non-negative while effectively capturing directional variations in diffusion. A value of 0 represents no diffusion, whereas 1 corresponds to maximum diffusion.}
    \label{fig:Fast_PINNs_NR_flow_velocity_invariant}
\end{figure}

A square domain with $L_x = L_y = 1$ is chosen. At the left boundary, a unit pressure is prescribed while the middle third right boundary is prescribed with zero pressure. The remaining boundaries are prescribed with zero velocity as shown in Fig. ~\ref{fig:Fast_PINNs_flow_problem_setup_pressure} A). The body force $\mathbf{b(x)}$ is assumed to be zero. The coefficient of permeability is taken as $k = 1$ over the domain. Flow subproblem is solved utilizing the governing equations in Eqs. ~\eqref{Eqn:PINNS_Fast_BoLM}--\eqref{Eqn:PINNS_vBC} solving for the pressure profile as presented in the Fig. ~\ref{fig:Fast_PINNs_flow_problem_setup_pressure} B) with the maximum pressure at the left boundary and minimum at the middle third of right boundary, satisfying the boundary conditions. 
%------------------------------------------;
%  Figure 10: Flow velocity problem setup  ;
%------------------------------------------;
The implemented PINNs code also computes the velocity field across the domain. For the given problem setup, the flow remains uniform, moving from the high-pressure region at the left boundary toward the low-pressure region at the middle third of the right boundary, as depicted in Fig.~\ref{fig:Fast_PINNs_NR_flow_velocity_invariant} A). As the flow approaches the top and bottom right boundaries, the velocity gradually diminishes and redirects toward the middle due to the imposed zero-velocity boundary conditions.

\begin{figure}[hbt!]
    \centering
    \includegraphics[scale=0.65]{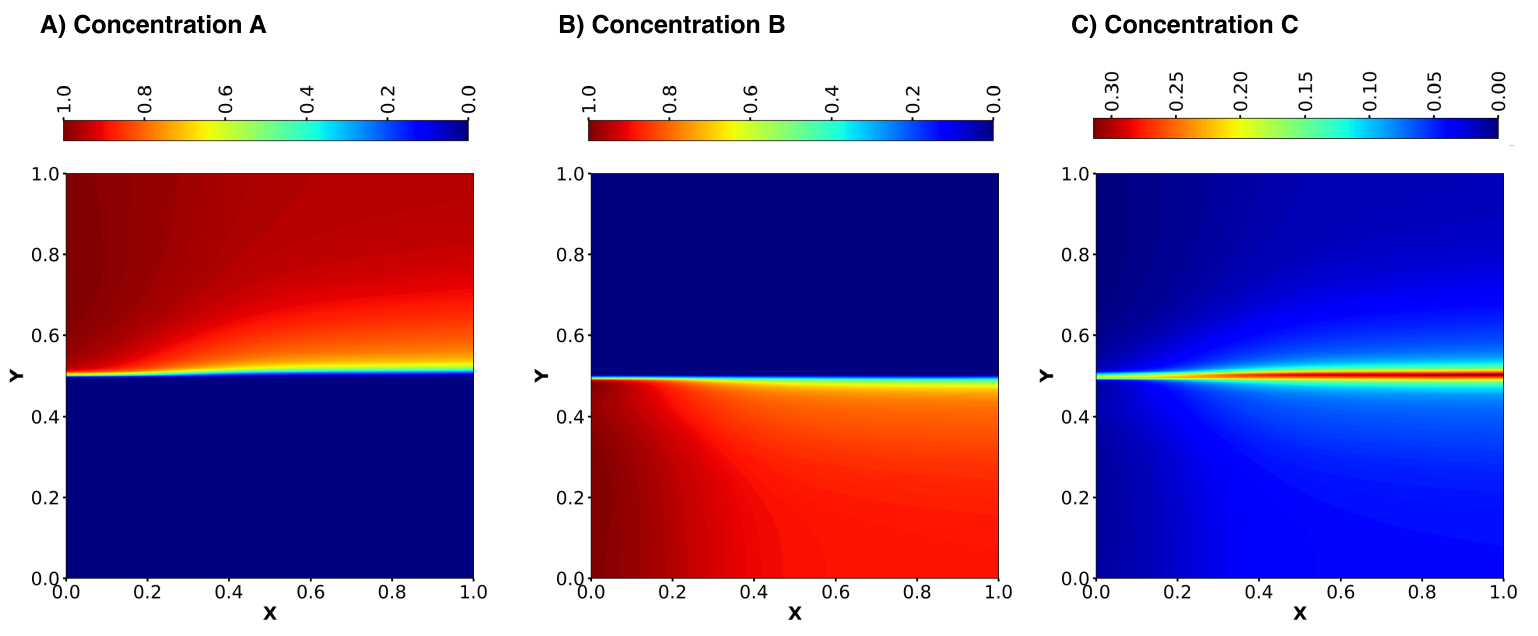}
    \caption{\textsf{Chemicals concentration:} The figure displays the concentration profiles of reactants A, B, and product C in (A), (B), and (C), respectively. The concentration of A is predominant in the top half, while B is concentrated in the bottom half, reflecting their respective introduction in the boundary value problem. The concentration values range from 0 to 1, with the highest concentration observed at the left boundary. The concentration of product C clearly illustrates its formation at the center, where reactants A and B interact. It begins forming at the left boundary and becomes more pronounced as it progresses toward the right boundary. The concentrations of all chemicals remain non-negative, adhering to the maximum principle.}
    \label{fig:Fast_PINNs_NR_flow_velocity_concentraton}
\end{figure}

The velocity field obtained from the flow subproblem is used to compute the mechanical dispersion tensor using Eq. ~\eqref{Fast_PINNs_diffusivity_tensor}. Two uncoupled equations in terms of the invariants $\Psi_A$ and $\Psi_B$ represented in Eqs.~\eqref{Fast_PINNs_invariant1} and ~\eqref{Fast_PINNs_invariant2} are solved using PINNs and plotted in Fig.~\ref{fig:Fast_PINNs_NR_flow_velocity_invariant} B) and C) respectively. The values of invariants $\Psi_A$ and $\Psi_B$  range from 0 to 1. $\Psi_A$ corresponds to the concentration of reactant A and product C, exhibiting maximum diffusivity in the upper half of the domain where reactant A is dominant, as per the problem setup. Likewise, $\Psi_B$ represents reactant B and product C, showing peak diffusion in the lower half, where reactant B is introduced. The invariants' plots provide insights into the diffusion behavior of different chemicals and their likely concentration distributions. 
%----------------------------------------;
%  Figure 11: Flow velcity and invariant ;
%----------------------------------------;
The concentrations of reactants A and B, along with product C, are determined using the invariants $\Psi_A$ and $\Psi_B$ based on the equations in Eqs.~\eqref{Fast_PINNs_cA}, ~\eqref{Fast_PINNs_cB}, and ~\eqref{Fast_PINNs_cC}. These concentrations are depicted in Fig.~\ref{fig:Fast_PINNs_NR_flow_velocity_concentraton}. Reactant A is most prominent in the upper half, where it was introduced in the problem setup, while reactant B is concentrated in the lower half, as it was introduced at the lower left boundary. The concentration values for both chemicals range from 0 to 1, with 1 representing the highest concentration at the left boundary and 0 indicating absence. Product C begins forming at the intersection of reactants A and B, emerging from the left boundary and becoming more prominent toward the right boundary. For a fast bimolecular reaction with uniform flow velocity, PINNs effectively capture the concentration distributions of both reactants and the product.
%-------------------------------------;
%  Figure 12: chemical concentration  ;
%-------------------------------------;
\subsection{Non-uniform explicit velocity}
\begin{figure}[hbt!]
    \centering
    \includegraphics[scale=0.65]{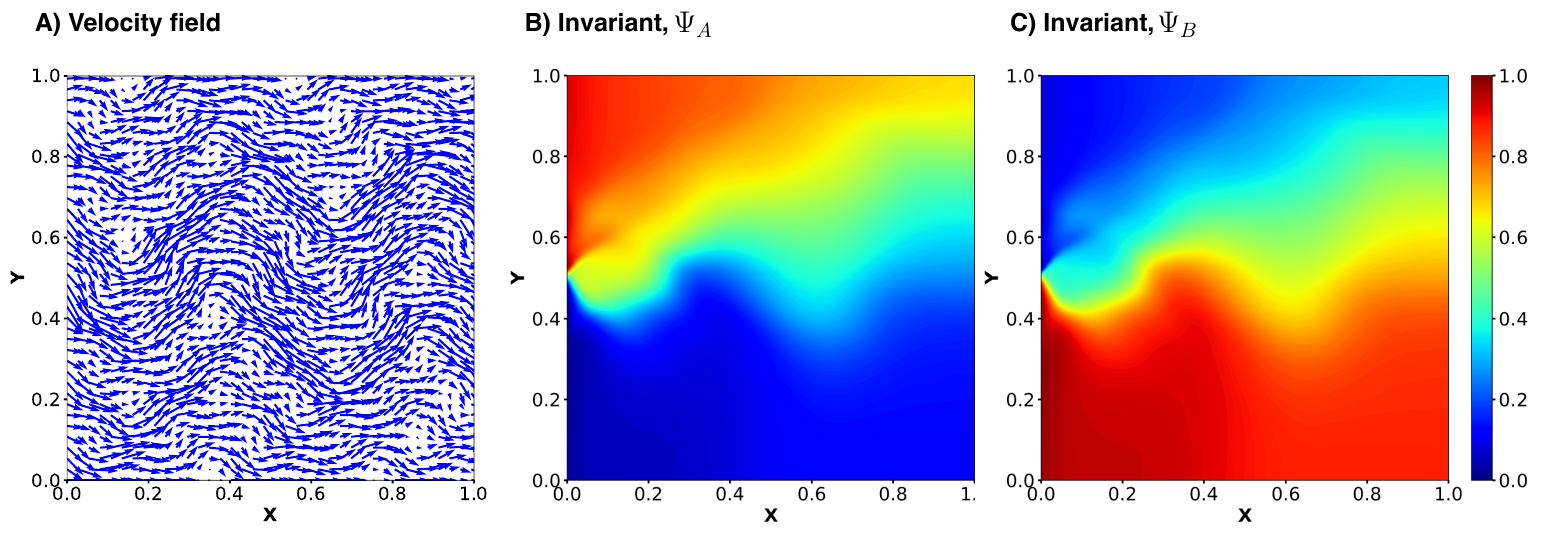}
    \caption{\textsf{Explicit velocity field and invariants: }The figure presents the solenoidal velocity field in (A), which is non-uniform with varying magnitudes and random directions. The diffusivity invariants are shown in (B) and (C), with values ranging from 0 to 1, ensuring compliance with the maximum principle. Higher values indicate increased diffusion in a particular direction, while lower values represent reduced or negligible diffusion.}
    \label{fig:Fast_PINNs_explicit_velocity_field_invariants}
\end{figure}

After evaluating the performance of PINNs on a fast bimolecular reaction with uniform flow velocity, the study is extended to examine the case of random and non-uniform flow velocity. To achieve this, an explicit velocity is explicitly prescribed and is selected referring from \citet{mudunuru2016enforcing}, with its $x$ and $y$ components, $\mathbf{v}(x)$ and $\mathbf{v}(y)$, defined as:
\begin{align}
    \label{explicit_velocity_x}
    & \mathbf{v}(x) = 1 + \sum_{i=1}^{3} A_i \, \frac{q_i \pi}{L_y} \, \cos \Big(\frac{p_i\pi \mathrm{x}}{L_x} - \frac{\pi}{2} \Big) \, \cos \Big( \frac{q_i \pi \mathrm{y}}{L_y}\Big) \\
    \label{explicit_velocity_y}
    & \mathbf{v}(y) = \sum_{i=1}^{3} A_i \, \frac{p_i \pi}{L_x} \, \sin \Big(\frac{p_i\pi\mathrm{x}}{L_x} - \frac{\pi}{2} \Big) \, \sin \Big( \frac{q_i \pi \mathrm{y}}{L_y}\Big)
\end{align}
where $(x,\, y)$ are spatial coordinates in the computational domain, $(p_1,\,p_2,\,p_3) = (4,5,10)$, $(q_1,\,q_2,\,q_3) = (1,5,10)$, and $(A_1,\,A_2,\,A_3) = (0.08,0.02,0.01)$. The velocity field is plotted with the arrow representing the velocity vector at the each collocation points in Fig.~\ref{fig:Fast_PINNs_explicit_velocity_field_invariants} A). The direction and magnitude of arrows clearly depicts the randomness in the velocity.  

The velocity field obtained is used to compute the diffusivity tensor using Eq. ~\eqref{Fast_PINNs_diffusivity_tensor}. The invariants $\Psi_A$ and $\Psi_B$ represented in Eqs.~\eqref{Fast_PINNs_invariant1} and ~\eqref{Fast_PINNs_invariant2} are solved using PINNs and plotted in Fig.~\ref{fig:Fast_PINNs_explicit_velocity_field_invariants} B) and C) respectively. Similar to the flow velocity case, the values of invariants $\Psi_A$ and $\Psi_B$  range from 0 to 1. $\Psi_A$ corresponds to the concentration of reactant A and product C, exhibiting concentration in the upper half of the domain where reactant A is dominant, as per the problem setup. Similarly, $\Psi_B$ represents reactant B and product C, exhibiting concentration in the lower half, where reactant B is introduced in the problem setup. 
%--------------------------------------------;
%  Figure 13: Explicit velocity & invariants ;
%--------------------------------------------;

\begin{figure}[hbt!]
    \centering
    \includegraphics[scale=0.65]{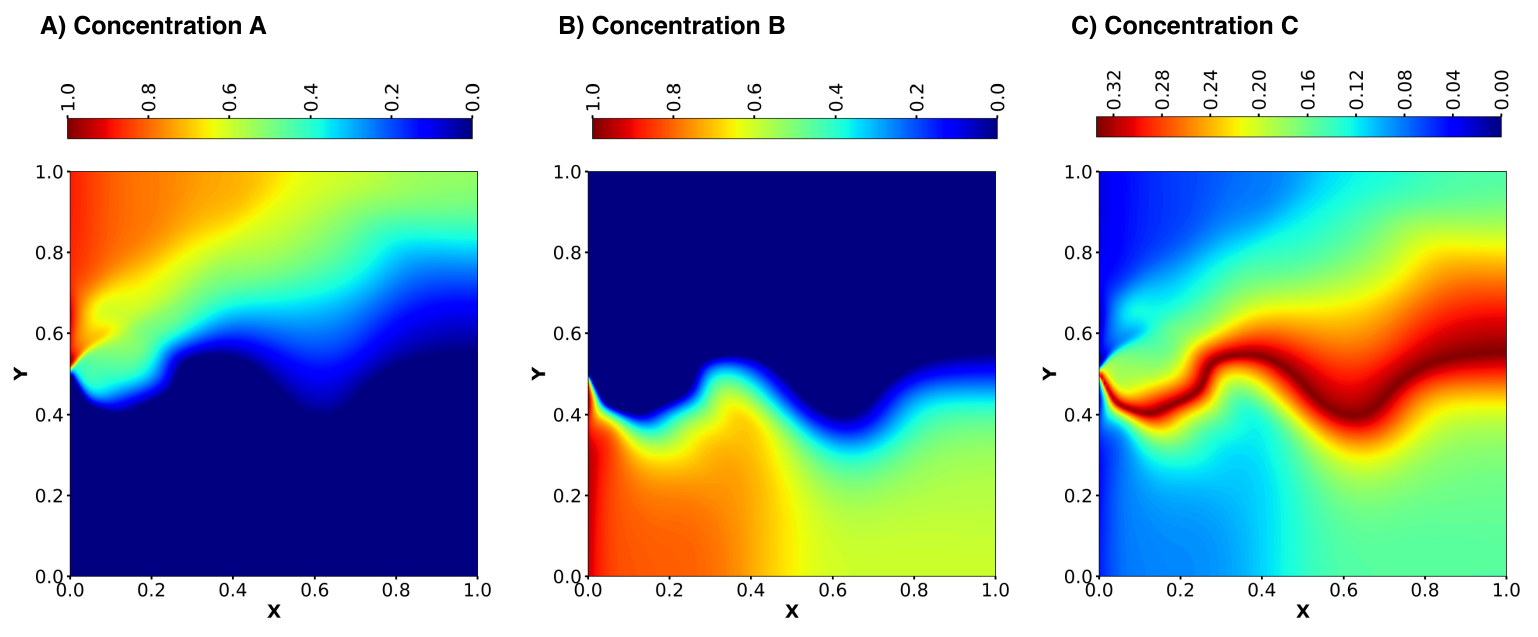}
    \caption{\textsf{Chemicals concentration and plume formation: }The figure presents the concentration profiles of reactants A, B, and product C in (A), (B), and (C), respectively. Reactant A is primarily concentrated in the top half, while B dominates the bottom half, consistent with their prescribed boundary conditions. The concentration values range from 0 to 1, adhering to the maximum principle, with the highest concentration observed at the left boundary. The concentration of product C distinctly illustrates its formation at the center, where reactants A and B interact. It initially emerges at the left boundary and gradually develops into a striking plume as it extends toward the right boundary.}
    \label{fig:Fast_PINNs_explicit_velocity_chemical_concentration}
\end{figure}

Finally, the concentration of reactants A and B along with product C are computed using the invariants $\Psi_A$ and $\Psi_B$ and presented in Fig.~\ref{fig:Fast_PINNs_explicit_velocity_chemical_concentration} . The mathematical calculation of the concentration is done using Eqs.~\eqref{Fast_PINNs_cA}--\eqref{Fast_PINNs_cC}. Reactant A exhibits its concentration at the top left boundary, gradually decreasing toward the right boundary. Similarly, reactant B is concentrated at the bottom left boundary and diminishes as it moves rightward. The formation of product C is observed at the interaction zone of reactants A and B, displaying a striking wavy concentration pattern. Product C begins forming at the left boundary and develops into a distinct plume extending toward the right boundary. Notably, the plume shifts slightly upward within the domain, which aligns with expectations, given that the coefficient $n_B$ is twice that of $n_A$, facilitating the upward movement of the plume.
%-------------------------------------;
%  Figure 14: Chemical concentration  ;
%-------------------------------------;

Despite the presence of a non-uniform velocity field, PINNs successfully capture the fast bimolecular reaction dynamics, revealing complex and visually striking nonlinear concentration distributions, comparable to those reported in \citet{mudunuru2016enforcing} with a similar problem setup. This further validates the robustness of the developed PINNs framework in solving chemical reaction problems under both uniform and non-uniform velocity conditions, answering the third science question. This makes PINNs well-suited for applications involving randomized velocity fields, such as subsurface flow.

%*********************************************;
%                                             ;
%  NAME                                       ;
%    S5_PINNS_Fast_Closure.tex                ;
%                                             ;
%*********************************************;
\section{CLOSURE}
\label{Sec:S6_PINNS_Fast_Closure}
We have developed a physics-informed modeling framework for fast bimolecular reactions, incorporating tests for flow and reactive transport in porous media. In applications where fluid flow uncertainties and the risk of chemical contamination are concerns, PINNs—a novel machine learning framework—offers a powerful solution. Unlike traditional data-driven methods, PINNs is trained using the governing physics rather than relying on large datasets. For fast bimolecular reactions, traditional FEM and FVM solvers typically require fine spatial discretization, stabilization, or specialized treatments to resolve sharp mixing-limited fronts. In contrast, the proposed PINNs framework offers a mesh-free and flexible alternative that enforces governing equations and reaction constraints directly at collocation points. While not intended to outperform classical solvers in large-scale production settings, PINNs are particularly attractive for data-scarce problems, rapid prototyping, inverse modeling, and scenarios where repeated mesh generation and solver tuning become prohibitive.

This approach has broad applicability, ranging from critical minerals to subsurface energy systems such as geothermal energy and geological carbon storage, where fluid flow and chemical transport play critical roles. A key challenge in such subsurface applications is the scarcity of sufficient data. PINNs address this limitation by leveraging the system of PDEs and boundary conditions, enabling accurate modeling without the need for extensive experimental or field data. PINNs is trained by minimizing the total loss, which is a combination of individual losses calculated based on how well the solution at each collocation point satisfies the given PDE and boundary conditions. The system of equations varies depending on the nature of the problem—for instance, whether it involves flow or transport, as in this study. The overall computation time is problem-dependent, with the simplest configurations completing in approximately 2 minutes and the most refined cases in under an hour. During the development of the PINNs framework and solving the problems, several key observations were made, which are summarized below:
\begin{enumerate}[(O1)]
    \item Accurate implementation of the PDEs and boundary conditions is crucial in the code. Even a minor error in formulating the system of equations can result in training the model on incorrect physics, ultimately leading to inaccurate solutions.
    \item When training the PINNs framework, particularly with a focus on PDEs and boundary conditions, the weighting of individual loss terms plays a crucial role. Assigning equal weights to all losses may cause the model to prioritize satisfying one while neglecting others. Therefore, the weight distribution should be carefully adjusted based on the relative magnitude of each loss to ensure that the final solution effectively satisfies the entire system of equations.
    \item Unlike conventional data-driven machine learning, PINNs require a balanced approach to loss minimization. Aggressively reducing the loss to an extremely small value can cause the model to focus solely on optimization while neglecting the underlying physics, potentially leading to solution divergence. Therefore, it is essential to maintain a careful balance between minimizing the loss and ensuring that the solution remains consistent with the governing physics.
\end{enumerate}
After successfully developing a functional PINNs framework, the model was tested on various phenomena, particularly flow and transport, to assess its solving capability and accuracy. It was then applied to model the fast bimolecular reaction to solve for chemical concentration. The key conclusions from this study are listed below:
\begin{enumerate}[(C1)]
    \item Three different patch tests—vertical, horizontal, and inclined—were conducted using a mixed formulation to model flow through porous media with heterogeneous permeability. The PINNs framework successfully predicted pressure and velocity, with results validated against the FEM solution. However, determining the optimal loss weight combination posed a challenge when solving with PINNs, making FEM the more straightforward approach in this case. Nonetheless, the PINNs framework maintained high accuracy, which remains a crucial aspect of modeling.
    \item The transport problem was tested using the PINNs framework on a square domain with a hole. PINNs successfully solved for the concentration while adhering to the governing PDE and boundary conditions. The results were compared with those obtained from FEM. A well-known drawback of FEM in diffusion problems is the violation of the non-negativity constraint (maximum principle). However, PINNs not only provided accurate solutions but also inherently preserved the maximum principle without any additional enforcement. This presents a significant advantage for diffusion-type problems, demonstrating the strength of PINNs in this test.
    \item A comprehensive uniform flow-transport-reaction problem was solved using the PINNs framework. The model effectively captured the fast bimolecular reaction and accurately predicted the concentrations of reactants and the resulting product. The obtained results closely align with those from a previous study \citep{mudunuru2016enforcing}, which used FEM to solve the fast bimolecular reaction. This demonstrates that PINNs is capable of solving an entire sequence of processes, from fluid flow to chemical reactions.
    \item Finally, a test problem with a random advection velocity field was solved using PINNs to evaluate its ability to handle a fast bimolecular reaction with a non-uniform flow. PINNs successfully predicted the concentrations of the reactant chemicals and product across the domain, resulting in a striking plume formation of the product. The chemical concentration plots obtained closely resemble those reported in \citet{mudunuru2016enforcing}. This further reinforces the capability of the novel PINNs framework in solving problems involving non-uniform flows with random velocities.
\end{enumerate}

In summary, PINNs excelled in all tests, successfully solving the fast bimolecular reaction, regardless of whether the flow velocity was uniform or random. This framework has the potential to be a game-changing tool for a wide range of applications involving critical minerals and materials. It can be applied across various fields, from surface water to subsurface applications, from conventional to unconventional energy solutions. While this study focuses on 2D domains for validation purposes, the PINNs framework is inherently extendable to 3D. In addition, a worthwhile extension of this work would be developing a complex reactive transport benchmark problem to represent wider applications such as the MoMaS \citep{carrayrou2010reactive}. Exploring this extension, along with comparative performance analysis against traditional reactive transport solvers \citep{de2010global, carrayrou2010looking, lagneau2010hytec,ahusborde2019fully}, presents a promising direction for future work aimed at evaluating scalability and practical deployment. 

%================================;
%  Section: Computational tools  ;
%================================;
\section*{COMPUTATIONAL TOOLS AND IMPLEMENTATION ENVIRONMENT}
For PINNs code development, we used \textbf{Spyder} (Anaconda Navigator 2.6.6; \url{https://www.anaconda.com/products/navigator}) as code editors. The implementation was carried out in \textbf{Python} (version 3.11; \url{https://www.python.org}) using the \textbf{JAX} library (version 0.4.28; \url{https://docs.jax.dev/en/latest}) for efficient numerical computation. For visualization, we used \textbf{Matplotlib} (version 3.7.2; \url{https://matplotlib.org}). Finite element simulations were performed using \textbf{COMSOL Multiphysics} (version 5.0; \url{https://www.comsol.com}), utilizing its built-in input/output features. Conceptual drawings were created using \textbf{Affinity Designer} (version 2; \url{https://affinity.serif.com}).

% %======================================;
% %  Appendix (Supplementary materials)  ;
% %======================================;
% \section*{\textbf{SUPPLEMENTARY MATERIAL}}
% %%
% Supplementary material and associated figures can be found in a separate file.

%=====================;
%  Data availability  ;
%=====================;
\section*{DATA AND CODE AVAILABILITY}
Data supporting the findings of this study are available from the corresponding authors (KBN) on a reasonable request.
The codes used in this paper will be made available and hosted at EMSL-Computing GitHub (\url{https://github.com/EMSL-Computing/AI4Minerals}) upon publication.

\vspace{0.1in}
%=========================;
%  Authors contributions  ;
%=========================;
\noindent\textbf{Author Contribution Declaration.} KBN conceived the study and developed the conceptual framework. Methodology was designed collectively by all authors (KA, MM, MKM, and KBN). KA conducted the formal analysis and performed the investigations with support from KBN. Resources were provided by KBN, MM, and MKM. KA and KBN prepared the original draft of the manuscript, and all authors (KA, MM, MKM, and KBN) contributed to the review and editing. Visualization was handled by KA. KBN supervised the project, managed its administration, and acquired funding. All authors have read and approved the final manuscript.

\vspace{0.1in}
%=======================;
%  Funding declaration  ;
%=======================;
\noindent\textbf{Funding Declaration.} The authors acknowledge the support from the Environmental Molecular Sciences Laboratory (EMSL), a DOE Office of Science User Facility sponsored by the Biological and Environmental Research program under contract no: DE-AC05-76RL01830 (Large-Scale Research User Project No: 60720, Award DOI: \url{10.46936/lser.proj.2023.60720/60008914}). KA and KBN also acknowledge the support from the UH Additive Manufacturing Institute (AMI) and the UH-Chevron Energy Fellowship.

\vspace{0.1in}
%==============;
%  Disclaimer  ;
%==============;
\noindent\textbf{Disclaimer}
This research work was prepared as an account of work sponsored by an agency of the United States Government. Neither the United States Government nor any agency thereof, nor any of their employees, makes any warranty, express or implied, or assumes any legal liability or responsibility for the accuracy, completeness, or usefulness of any information, apparatus, product, or process disclosed, or represents that its use would not infringe privately owned rights. Reference herein to any specific commercial product, process, or service by trade name, trademark, manufacturer, or otherwise does not necessarily constitute or imply its endorsement, recommendation, or favoring by the United States Government or any agency thereof. The views and opinions of authors expressed herein do not necessarily state or reflect those of the United States Government or any agency thereof.

%================;
%  Bibliography  ;
%================;
\bibliographystyle{plainnat}
\bibliography{Master_References}
\end{document}